\documentclass[]{elsarticle} 

\usepackage[utf8]{inputenc} 

\usepackage{geometry} 
\RequirePackage{amsfonts,amssymb,stmaryrd}
\RequirePackage{mathtools, xcolor}
\usepackage{breqn, siunitx}
\RequirePackage{bm, graphicx, xspace}
\usepackage{commands}
\usepackage{ulem}

\usepackage{booktabs} 
\usepackage{array} 
\usepackage{paralist} 
\usepackage{verbatim} 
\usepackage{subfig} 

\usepackage{fancyhdr} 
\usepackage{sectsty}
\allsectionsfont{\sffamily\mdseries\upshape} 

\usepackage[nottoc,notlof,notlot]{tocbibind} 
\usepackage[titles,subfigure]{tocloft} 

 \graphicspath{{figures/}}
 \usepackage[hyphens]{url}
 
\title{A fully coupled numerical model of thermo-hydro-mechanical processes and fracture contact mechanics in porous media}

\author[1]{Ivar Stefansson\corref{cor1}
	}
\ead{ivar.stefansson@uib.no}
\author[1]{Inga Berre} 
\ead{inga.berre@uib.no}

\author[1]{Eirik Keilegavlen}
\ead{eirik.keilegavlen@uib.no}

\cortext[cor1]{Corresponding author}
\address[1]{University of Bergen, Bergen,
Norway}
\date{} 

\begin{document}

\begin{abstract}
Various phenomena in the subsurface are characterised by the interplay between \is{deforming structures such as fractures and} coupled thermal, hydraulic and mechanical processes. Simulation of subsurface dynamics can provide valuable phenomenological understanding, but requires models which faithfully represent the dynamics involved; these models therefore are themselves highly complex.

This paper presents a mixed-dimensional thermo-hydro-mechanical model designed to capture the process--structure interplay using a discrete--fracture--matrix framework. It incorporates tightly coupled thermo--hydro--mechanical processes based on \is{balance} laws for momentum, mass and \is{energy} in subdomains representing the matrix and the lower-dimensional fractures and fracture intersections.
The deformation of explicitly represented fractures is modelled by contact mechanics relations and a Coulomb friction law, with a novel formulation consistently integrating fracture dilation in the governing equations. 

The model is discretised using multi-point finite volumes for the balance equations and a semismooth Newton scheme for the contact conditions and is implemented in the open-source fracture simulation toolbox PorePy. Finally,  simulation studies demonstrate the model's convergence, investigate process--structure coupling effects, explore different fracture dilation models and show an application of the model to stimulation and long-term cooling of a three-dimensional geothermal reservoir.
\end{abstract}

\maketitle

\section{Introduction}
Fluid injection operations into the subsurface are common in e.g.~geothermal energy and petroleum production, wastewater disposal, CO$_2$ storage and groundwater management. Injection can severely alter subsurface hydraulic, mechanical, thermal and chemical conditions. These coupled processes are strongly affected by preexisting fractures, which represent extreme heterogeneities and discontinuities in the formation. The processes may in turn cause deformation of the fractures, giving rise to dynamic and highly complex process--structure interactions.  

 In some subsurface engineering operations, fracture deformation is deliberately induced, e.g.~to enhance permeability through hydraulic stimulation, in which fluid is injected at elevated pressure to overcome a fracture's frictional resistance to slip \cite{pine1984downward,rutqvist2003role,evans2005microseismicity}. There may also be interest in preventing deformation of fractures to, for example, avoid  induced seismicity of unacceptable magnitude in disposal of wastewater  \cite{keranen2014sharp,dewaal2015production,improta2015detailed, keranen2018induced} or during hydraulic stimulation of fractured geothermal reservoirs \cite{dorbath2009seismic,ellsworth2019triggering, haring2008characterisation}.

As data related to subsurface dynamics are limited, physics-based modelling can complement data analysis in understanding governing mechanisms for fracture deformation. This requires numerical simulation tools that can capture the governing structure of the fractured formation and relevant coupled processes as well as process--structure interactions, which necessitates explicit representation of both the matrix and dominant fractures in the model. Typically, major fractures or faults are represented explicitly while the rest of the domain is represented as a matrix continuum, possibly integrating effects of finer-scale fractures.

In a spatial grid, there are two alternatives for representing such a discrete--fracture--matrix (DFM) model: Resolving the width of the fractures in the grid in an equi-dimensional model imposes severe restrictions put on the spatial discretisation of the domain due to the high aspect ratio of the fractures, thereby limiting the number of fractures that can be included in the model. A geometrically simpler alternative, which was introduced for flow models, is a co-dimension one model, where fractures are represented as objects of one dimension lower than the surrounding domain  \cite{Karimi-Fard2003efficient, Martin2005, reichenberger2006mixed, nordbotten2019unified}. In contrast to simulation models for coupled flow and mechanics that treat faults as equidimensional zones of different rheology resolved in the grid \cite{rutqvist2008coupled, rutqvist2013modeling, wassing2014coupled}, the co-dimension DFM model facilitates modelling of fracture slip and dilation \cite{cappa2011modeling, ucar2018three} and can be combined with full mechanical fracture opening \cite{ferronato2008subsidence, gallyamov2018discrete}. A conceptually simpler alternative to co-dimension one DFM models is to incorporate only the dynamics in the fracture network and either disregard the dynamics in the matrix altogether or approximate them using semi-analytical methods. These approaches are based on Discrete-Fracture-Network (DFN) representations \cite{shapiro1985simulation, long1985model} and will be referred to as DFN methods.

Driven by the need to improve the result of injection operations and avoid unacceptable environmental impacts, intense focus has been placed on physics-based modelling. Early works by \citet{willis1996progress, rahman2002shear_dilation, kohl2007predictive} and \citet{bruel2007using} developed  DFN-type models considering only deformation and flow in the fractures and using a Coulomb friction law to model fracture slip due to changes in effective stress as a consequence of local change in fluid pressure. Later, \citet{baisch2010numerical} improved on this type of model by including redistribution of shear stress along the fracture as a consequence of slip through a block-spring model.  \citet{mcclure2011investigation} further developed the modelling of mechanical interaction between fractures with the boundary integral equation method and introduced a rate-and-state friction model. This type of method has been combined with fracture propagation \cite{norbeck2016embedded, ciardo2019effect}. As only the fracture is discretised when using the boundary integral equation method, models based on this approach can be classified as DFN-type models. Common to all of these approaches is use of semi-analytical approaches and sequential coupling of physical processes.    

The last decade has seen developments in the inclusion of dynamics in the matrix as well as improved models and numerical solution schemes for coupling of different dynamics. Building on previously developed DFN-type models, \citet{mcclure2011investigation} and \citet{mcclure2015generation} introduced a semi-analytical leakoff term to mimic fracture--matrix flow. \citet{norbeck2016embedded} expanded on previous models developed by \citet{mcclure2011investigation} and accounted for the interaction between fracture and matrix flow through an embedded discrete fracture model, where flow in the fracture and in the matrix are discretised on non-conforming grids and connected through transfer terms. 
Hydro-mechanical simulation tools based on co-dimension one DFM models combined with Coulomb friction laws for fracture slip  have also been introduced \cite{ucar2018three, franceschini2020stabilisedkkt, berge2020hmdfm, jha2014coupled, garipov2016discrete, settgast2017fracking, ucar2017postinjection}, motivated by applications related to CO$_2$-storage \cite{jha2014coupled}, oil and gas production \cite{garipov2016discrete, settgast2017fracking} and hydraulic stimulation of fractured geothermal reservoirs \cite{ucar2017postinjection, ucar2018three}. In these tools, the contact mechanics conditions are typically handled through Lagrange-multipliers  \cite{ jha2014coupled, berge2020hmdfm, franceschini2020stabilisedkkt} or penalty methods \cite{garipov2016discrete}; see e.g.~\citet{wriggers2004computational}.

More recently, thermal effects have been taken into account in deformation of fractured porous media. Based on a DFN-type model, where the boundary integral equation method was used so that only the fracture is discretised, \citet{ghassemi2011three} included thermo-poroelastic effects in the matrix. Based on a DFM conceptual model, \citet{pandey2017coupled} and \citet{salimzadeh2018thm} have presented models with linear thermo-poroelasticity for the matrix combined with flow, heat transfer and deformation of a single fracture. However, none of these works included modelling of fracture slip or shear dilation when fracture surfaces are in contact. \citet{gallyamov2018discrete} consider a conceptually similar model which includes multiphase flow and a fracture-contact-mechanics model combined with opening and propagation of fractures, and present simulation studies with a large number of fractures. Their work considers the impact of the contact traction on the hydraulic aperture of closed fractures. In contrast to the majority of previously mentioned works, where simplifications that impact the solution are made in the solution of the coupled system, they solve the equations fully coupled, i.e.\ the flow, energy and mechanics equations are solved simultaneously building on the work by \citet{garipov2016discrete}. 

Recent work by \citet{garipov2019nitschewaterflooding} combines several previous developments. Their work is based on a DFM model and considers a fully coupled thermo-poroelastic model for the matrix, flow and heat transfer in the fractures and contact mechanics for fractures based on a Coulomb friction law. Energy and mass conservation are discretised by a finite volume  method, while momentum is discretised by a Galerkin finite element method. This work also presents robust treatment of couplings in the model. However, while the work accounts for permeability enhancement due to full opening of fractures as well as shear dilation, stress response due to dilation as a consequence of slip is not included in the model.

This paper presents a mathematical model based on a mixed-dimensional DFM representation of coupled thermo-hydro-mechanical (THM) processes in a porous rock containing deforming fractures with an accompanying discretisation and numerical solution approach. 
The model fully couples fluid flow and transport in both matrix and fractures, linear thermo-poromechanics in the matrix and nonlinear fracture deformation. Fracture deformation is based on traction balance, nonpenetration and a Coulomb type friction law, and allows for shear slip and dilation as well as complete fracture opening. To the authors' knowledge, this is the first model that consistently and fully coupled represents stress redistribution due to slip-induced dilation of fractures.  As demonstrated by the numerical results, the effect of this coupling can be significant. 

Based on the modelling of fractures as lower-dimensional surfaces, the domain is decomposed into subdomains of different dimensions corresponding to matrix, fractures and intersections. Model equations, sets of variables and parameters
are defined on each subdomain and the interfaces between them. The resulting mixed-dimensional model \cite{keilegavlen2020porepy} 
facilitates systematic modelling on the decomposed structure while incorporating interaction between processes both within and between subdomains. 
The governing balance equations in each subdomain are discretised based on multi-point finite volume methods preserving local conservation, using the same spatial grid for discretisation of all processes. The contact conditions are formulated in terms of contact tractions at the interfaces between fractures and matrix blocks, similar to approaches based on Lagrange-multipliers \cite{wriggers2004computational}. The nonlinear fracture deformation equations are discretised using a semismooth Newton scheme formulated as an active set method. 
The model and its implementation extend the work presented by \citet{berge2020hmdfm}, who \is{consider the poromechanical problem without flow in the fractures,} shear dilation and thermal effects.

The model is presented in Section \ref{sec:model}, and its discretisation is described in Section \ref{sec:discretisation}. In both sections, particular emphasis is placed on fracture deformation as well as its impact on the balance equations for the fractures and the back-coupling to the higher-dimensional momentum balance. 
Three examples are presented in Section \ref{sec:results}: The first investigates governing mechanisms and coupling effects and verifies the model and its implementation in a convergence study. In the second, three different models for fracture dilation are compared. In the last example, the model is applied to a \is{three-dimensional} hydraulic stimulation and long-term cooling scenario for geothermal energy extraction.
Finally, Section \ref{sec:conclusion} provides some concluding remarks.

\section{\is{Mixed-dimensional} governing equations}\label{sec:model}
This section describes the model for THM processes in a porous medium with contact mechanics at the fractures.
It relies on a DFM model in which the matrix, the fractures and fracture intersections are explicitly
represented by individual subdomains. To avoid resolving the small geometric distances introduced to the fracture network
geometry by the high aspect ratio of the fractures, the dimensions of the fracture and intersection subdomains are reduced.
The subdomains are collected in a hierarchical structure and connected by interfaces to yield the full mixed-dimensional model. 

Decomposition into subdomains facilitates tailored modelling of processes in distinct subdomains,
while interactions between subdomains take place on the interfaces.
Specifically, separate sets of variables, equations and parameters
are defined on each subdomain and interface.
This procures the flexibility needed to model the highly complex system
arising from the coupled THM system posed in both matrix and fractures.

The model consists of balance equations for momentum, mass and \is{energy} and relations governing the fracture deformation posed on the subdomains.
These are supplemented by constitutive laws and equations for coupling over the interfaces. 
The equations are formulated in terms of the primary variables displacement, pressure, temperature and contact traction on the fractures.

\is{The governing equations for THM-processes in a mono-dimensional porous medium are introduced succinctly in Section \ref{sec:matrix_equations}}, followed by a more elaborate presentation of  the lower-dimensional scalar  equations for deforming fractures and intersections emphasising the effect of volume change  in Section \ref{sec:lower_dim_TH}. Section \ref{sec:fracture_deformation} describes the model for fracture deformation and its relation to volume change.

\subsection{Matrix thermo-poromechanics}\label{sec:matrix_equations}
We first consider the governing equations in the matrix domain consisting of a solid and a fluid phase.
For the remainder of this subsection, let \domain[] denote the matrix domain,  \normal[] the outer normal vector on the boundary \boundary{}{} and $\dA=\normal \dx$.
Neglecting inertia, the momentum balance equation reads 
\begin{gather}\label{eq:momentum_balance}
\begin{aligned}
\int_{\boundary{}{}} \stress\cdot \dA = \int_{\domain[]} \vectorSource[\displacement]\dx
\end{aligned} 
\end{gather}
with \vectorSource[\displacement] denoting body forces and the thermo-poroelastic stress tensor for infinitesimal deformation modelled as linearly elastic obeying an extended Hooke's law
\begin{gather}\label{eq:mechanics_stress} 
  \begin{aligned}
    \stress-\stress[0] = \frac{\stiffness}{2}:(\nabla\displacement + \nabla\displacement^T)- \biotAlpha \left(\pressure-\pressure[0]\right) \identity - \thermalExpansion[s]  \bulkModulus{s} \left(\temperature - \temperature[0]\right) \identity. 
  \end{aligned}
\end{gather} 
Here, \stiffness denotes the stiffness tensor, \biotAlpha  the Biot coefficient,  \thermalExpansion[s] the volumetric thermal expansion and \bulkModulus{s} the bulk modulus of the solid, while \displacement, \pressure, and \identity are displacement, pressure and identity matrix. 
The subscript 0 indicates the reference state of a variable.
The assumption of local thermal equilibrium between fluid and solid leads to a single temperature unknown \temperature and to the definition of effective  coefficients, which are computed as weighted sums \cite{mctigue1986thermoelastic}: 
\begin{gather}\label{eq:effective_parameters}
\begin{aligned}
(\cdot)_{\eff} = \porosity (\cdot)_s  + \left(1-\porosity\right)(\cdot)_f. 
\end{aligned} 
\end{gather}
Subscripts $s$ and $f$ indicate that all quantities in the parentheses are for the solid and fluid, respectively. 
Herein, the relations $ \frac{\stiffness}{2}(\nabla\displacement + \nabla\displacement^T) = \shearModulus(\nabla\displacement + \nabla\displacement^T) + \bulkModulus{s} \trace{\nabla\displacement}\identity$ and
$\vectorSource[\displacement]=\density[s]\gravityVector$ are used, with  \shearModulus denoting the shear modulus, $\trace{\cdot}$ the trace of a tensor, \density[s] the solid density and \gravityVector the gravitational acceleration vector.

We assume the density of the slightly compressible fluid to follow
\begin{gather}\label{eq:density_eos}
\begin{aligned}
\density[f] = \density[0] \exp\left[\frac{1}{\bulkModulus{f}} (\pressure-\pressure[0]) - \thermalExpansion[f] (\temperature-\temperature[0])\right],
\end{aligned}
\end{gather}
with \thermalExpansion[f] and \bulkModulus{f} denoting fluid thermal expansion coefficient and bulk modulus, which is the inverse of the compressibility. Balance of mass reads \cite{coussy2004poromechanics}
\begin{gather}\label{eq:mass_balance}
\begin{aligned}
\int_{\domain[]}
\left(\frac{\porosity}{\bulkModulus{f}}+\frac{\biotAlpha-\porosity}{\bulkModulus{s}}\right) \dd{\pressure}{\timet} + \biotAlpha \dd{(\divergence\displacement)}{\timet}-\thermalExpansion[\eff]\dd{\temperature}{\timet} dx+\int_{\boundary{}{}} \fluidFlux \cdot \dA= \int_{\domain[]}\source[\pressure] \dx,
\end{aligned}
\end{gather}
with porosity \porosity and effective thermal expansion \thermalExpansion[\eff].
Fluid flux relative to the solid is denoted by \fluidFlux and volume sources and sinks by \source[\pressure]. With \permeability denoting the permeability and \viscosity the viscosity, the flux is modelled according to Darcy's law:
\begin{gather}\label{eq:mass_darcy}
\begin{aligned}
 \fluidFlux = -\frac{\permeability}{\viscosity} (\nabla \pressure - \density[]\gravityVector).
\end{aligned}
\end{gather} 

Neglecting viscous dissipation on an assumption of small velocities~\cite{coussy2004poromechanics}, the energy balance equation is \cite{tong2010fully,cacace2017flexible, coussy2004poromechanics}
\begin{gather}\label{eq:heat_balance1}
\begin{aligned}
\int_{\domain[]}\dd{}{\timet}\left[\porosity U_f+(1-\porosity)U_s\right]+ \thermalExpansion[s] \bulkModulus{s} \temperature[0] \dd{(\divergence\displacement)}{\timet} \dx+\int_{\boundary{}{}} (\advectiveHeatFlux[]+\conductiveHeatFlux)\cdot \dA=\int_{\domain[]}\source[\temperature]\dx.
\end{aligned}
\end{gather}
Here, the internal energy of the solid is $\internalEnergy[s] = \density[s] \heatCapacity[s] \temperature$, where $\heatCapacity[s]$ is the specific heat capacity of the solid. Assuming a simplified low-enthalpy description of the fluid, we approximate the fluid internal energy as $\internalEnergy[f] = \density[f] h_f$, where the fluid enthalpy is approximated as $h_f = \heatCapacity[f] \temperature$ and  $\heatCapacity[f]$ is the specific heat capacity of the fluid. The thermoelastic dissipation term involving the displacement represents the effect of the solid elastic deformation on the temperature distribution.
The total heat flux is the sum of the advective flux 
\begin{gather}\label{eq:advectiveflux}
\begin{aligned}
\advectiveHeatFlux[]=\density[f] \heatCapacity[f] \temperature \fluidFlux 
\end{aligned}
\end{gather}
and the diffusive Fourier flux 
\begin{gather}\label{eq:conductiveflux}
\begin{aligned}
\conductiveHeatFlux=-\heatConductivity[\eff] \nabla \temperature. 
\end{aligned}
\end{gather}
In the computation of the effective heat conductivity \heatConductivity[\eff] by Eq. \eqref{eq:effective_parameters}, dispersion due to micro-scale tortuous flow in the porous medium is neglected. The source term is assumed to equal the internal energy of the fluid of the volume source and sink terms, i.e.\ $\source[\temperature]=\density[f]\heatCapacity[f]\temperature\source[\pressure]$.

The energy equation can thus be written
\begin{gather}\label{eq:heat_balance2}
\begin{aligned}
 \int_{\domain[]}\dd{}{\timet}\Big[(\density \heatCapacity)_{\eff} \temperature\Big] + \thermalExpansion[s] \bulkModulus{s} \temperature[0] \dd{(\divergence\displacement)}{\timet} \dx+\int_{\boundary{}{}}  (\density[f]\heatCapacity[f] \temperature \fluidFlux-\heatConductivity[\eff] \nabla \temperature) \cdot \dA = \int_{\domain[]}\source[\temperature]\dx.
\end{aligned}
\end{gather}
where $(\density \heatCapacity)_{\eff}$ denotes the effective heat capacity of the porous medium and is calculated by Eq.~\eqref{eq:effective_parameters}.
The first term of Eq.~\eqref{eq:heat_balance2} can be expanded, giving
\begin{gather}\label{eq:heat_balance_1term_2term}
\begin{aligned}
\dd{}{\timet}\Big[(\density \heatCapacity)_{\eff} \temperature\Big] = (\density \heatCapacity)_{\eff}\dd{\temperature}{\timet}+\temperature \dd{(\density \heatCapacity)_{\eff}}{\timet} = (\density \heatCapacity)_{\eff}\dd{\temperature}{\timet}+\temperature\left[\left(\frac{\density\heatCapacity}{\bulkModulus{}}\right)_\eff \dd{p}{\timet}-\left(\density \heatCapacity \thermalExpansion[]\right)_\eff \dd{\temperature}{\timet}\right],
\end{aligned}
\end{gather}
where effective parameters are calculated using Eq.~\eqref{eq:effective_parameters}.

\begin{figure}[tp]
\centering
\includegraphics[width=.95\textwidth]{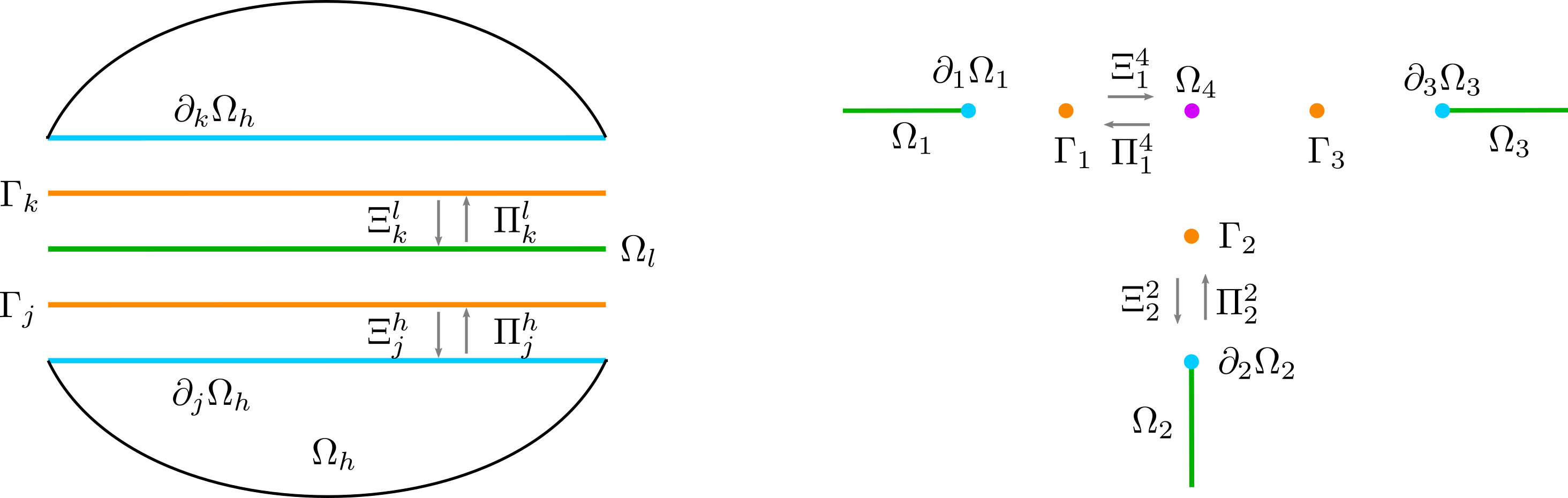}
\caption{Schematic representation of a throughgoing one-dimensional fracture in a two-dimensional matrix (left) and three fractures meeting at an intersection point (right).
All subdomains, internal boundaries and interfaces are indicated, as are select projection operators.
\domain[h] and  \domain[l] are separated by the interfaces $\interface[j]$ and $\interface[k]$ corresponding to the internal boundaries \boundary{h}{j} and \boundary{h}{k}. 
The intersecting fractures are indexed counterclockwise from 1 through 3 \is{and intersect} at \is{the }point \domain[4]. The interfaces are numbered so that \interface[i] matches \boundary{i}{i}. 
In the model, \domain[l], $\interface[j]$, $\interface[k]$,  \boundary{h}{j} and \boundary{h}{k} coincide geometrically, as do all zero-dimensional points in the right figure.}
\label{fig:fracture_schematic}
\end{figure}
\subsection{Lower-dimensional flow and heat transfer}\label{sec:lower_dim_TH}
This section derives balance equations  for mass and energy for fluid-filled fractures and intersections which may undergo significant relative deformation and volume change, giving rise to an additional term compared to equations for static domains.
Along with outlining dimension reduction for the mass and energy equations, the connection between the subdomains of the mixed-dimensional model is presented.

Some notation is needed for a unified description of the mixed-dimensional model of a fractured porous domain of dimension $\nd=3$ or $\nd=2$. The domain is split into connected subdomains corresponding to the
rock matrix, the co-dimension one fracture planes and co-dimension two fracture intersections. In the case $\nd=3$, the model also generalises to account for intersections of fracture intersection lines, i.e.\ zero-dimensional points. A subdomain is denoted by \domain[i] and its
boundary by \boundary{i}{}. The subscript $i$ is also used to identify variables defined within \domain[i], but suppressed as context allows.
Each part \boundary{i}{j} of the internal boundary is associated with an interface
\interface[j] to an immersed lower-dimensional domain \domain[l] (see \figRef{fracture_schematic}).
All lower- and higher-dimensional interfaces of a subdomain are collected in the sets 
\lowerSet[] and \higherSet[]; in particular, the interfaces corresponding to surfaces of fracture $i$ constitute \higherSet[i].
Where convenient, the higher- and lower-dimensional neighbours of an interface are denoted by \domain[h] and \domain[l], respectively. 

Finally, four types of projection operators are needed to transfer variables
between interfaces and the neighbouring higher- and lower-dimensional subdomains.
As illustrated in \figRef{fracture_schematic}, projection from the
interface to the subdomains is performed by \projectFromInterface{h}{j} and
\projectFromInterface{l}{j}, respectively, whereas \projectToInterface{h}{j} and
\projectToInterface{l}{j} project from the part of a subdomain geometrically coinciding with the interface to the interface. 

The thickness of a fracture is characterised by the aperture $\aperture$ [\si{\metre}], which will be related to the fracture deformation in Section \ref{sec:fracture_deformation}. The aperture of an intersection is taken to be the average of the intersecting higher-dimensional neighbours, i.e.
\begin{equation}\label{eq:aperture_inheritance}
\begin{aligned}
\aperture[l] = \frac{1}{|\higherSet[l]|}\sum_{j\in \higherSet[l]} \projectFromInterface{l}{j}\projectToInterface{h}{j} \aperture[h],
\end{aligned}
\end{equation}
where the projection operators transfer the higher-dimensional aperture first to the interface and then to the intersection, see Fig.~\ref{fig:fracture_schematic}.
The specific volume $\specificVolume[i]=\aperture^{\nRed}$ accounts for the dimension reduction from the deforming equi-dimensional \domainEq to the corresponding spatially fixed $d$-dimensional \domain[i]. 
With $\specificVolume=1$ for $d=\nd$, dimension reduction for a scalar quantity \scalar and a vector quantity \vect satisfy
\begin{gather}\label{eq:dimension_reduction_scalar}
\begin{aligned}
\dD{}{\timet}\int_{\domainEq}\scalar \dx &=\int_{\domain[i]}\dd{}{\timet} (\specificVolume[i] \scalar[i]) \dx
\end{aligned}
\end{gather}
and
\begin{gather}\label{eq:dimension_reduction_flux}
\begin{aligned}
\int_{\boundaryEq}\vect \cdot \dA &=\int_{\boundary{i}{}}\specificVolume[i]\vect[i] \cdot \dA
- \sum_{j\in \higherSet[i]} \projectFromInterface{l}{j} \bigg(\projectToInterface{h}{j} \specificVolume[h] \int_{\interface[j]} \iota_j \dx\bigg),
\end{aligned}
\end{gather}
respectively.
Here,  \vect[i] denotes the tangential $d$-dimensional flux, while the interface flux into the domain $\iota_j$ satisfies  $\projectFromInterface{h}{j} \iota_j = \vect[h] \cdot \normal[h]$ on \boundary{h}{j}, with \normal[h] denoting the outwards normal vector of \domain[h].
The weighting with \specificVolume[h] ensures that the interface flux matches the dimension of fluxes of the higher-dimensional neighbour, which are scaled by specific volumes as seen by the expression for the tangential flux.
Note that all differentials in reduced integrals should be interpreted as relative to the domain of integration; i.e.~\dx is two-dimensional for a fracture with $d=2$. 
Furthermore, $\dA=\normal[i]\dx$, where \normal[i] denotes the outwards normal at \boundary{i}{} lying in the tangent plane of \domain[i]. For $d=1$, the boundary integral equals evaluation of the integrand at the boundary points.

Assuming unitary fracture porosity, the fluid mass balance equation for a deforming equi-dimensional domain \domainEq is
\begin{gather}\label{eq:mass_balance_eq}
\begin{aligned}
\dD{}{\timet}\int_{\domainEq}\density[f] \dx
+\int_{\boundary{}{}}  \density[f]\fluidFlux \cdot \dA = \int_{\domainEq} \density[f]\source[\pressure] \dx.
\end{aligned}
\end{gather}
Averaging in the normal direction for the tangential flux and replacing the normal part of the boundary by \interface[j] according to Eq.~\eqref{eq:dimension_reduction_flux}, the fluid flux becomes
\begin{gather}\label{eq:fluid_flux_tangential_normal}
\begin{aligned}
 \int_{\boundaryEq}\fluidFlux\cdot\dA =\int_{\boundary{i}{}}\specificVolume[i]\fluidFlux[i]\cdot\dA - \sum_{j\in \higherSet[i]} \projectFromInterface{l}{j}\bigg(\projectToInterface{h}{j}\specificVolume[h]\int_{\interface[j]}\interfaceFluidFlux[j] \dx\bigg).
 \end{aligned}
\end{gather}
Thus, the interdimensional coupling between \domain[h] and \domain[l] takes the form of interface fluid fluxes \interfaceFluidFlux[j], which 
also appear as a Neumann condition for \domain[h]:
\begin{gather}\label{eq:interface_Neumann_condition}
    \begin{aligned}
		\fluidFlux[h] \cdot \normal[h]=\projectFromInterface{h}{j}\interfaceFluidFlux. 
	\end{aligned}
  \end{gather}
The interface flux is modelled using a Darcy type law extended from \citet{Martin2005} to account for gravity
  \begin{gather}\label{eq:interface_fluid_flux}
    \begin{aligned} 
\interfaceFluidFlux = - \frac{\permeability[j]}{\projectToInterface{l}{j}\viscosity[l]}\bigg( \frac{2}{\projectToInterface{l}{j}\aperture[l]} \left( \projectToInterface{l}{j} \pressure[l] - \projectToInterface{h}{j} \pressure[h]  \right) -\projectToInterface{l}{j}\density[f,l] \gravityVector \cdot \projectToInterface{h}{j}\normal[h] \bigg).
    \end{aligned}
  \end{gather}
  Both the weighting by \aperture[l] and \specificVolume[h] and the normal permeability \permeability[j] arise through dimension reduction. %
The remaining terms of Eq.~\eqref{eq:mass_balance_eq} are averaged in the normal direction using Eq.~\eqref{eq:dimension_reduction_scalar}, and Eq.~\eqref{eq:density_eos} is inserted for the fluid density. Collecting terms, dividing by \density[f] and assuming Darcy's law with tangential permeability \permeability yields the dimensionally reduced mass balance
\begin{gather}\label{eq:reduced_mass_balance_darcy}
\begin{aligned}
&\int_{\domain[i]} \specificVolume[i]\left(\frac{1}{\bulkModulus{f}}\dd{\pressure[i]}{\timet} - \thermalExpansion[f]\dd{\temperature[i]}{\timet}\right) +\dd{\specificVolume[i]}{\timet} \dx
-  \int_{\boundary{i}{}}\specificVolume[i]\frac{\permeability}{\viscosity}\left(\nabla\pressure[i]-\density[f]\gravityVector\right)\cdot\dA \\
&- \sum_{j \in \higherSet[i]}\projectFromInterface{i}{j}\bigg(\projectToInterface{h}{j}\specificVolume[h]\int_{\interface[j]}\interfaceFluidFlux[j]\dx \bigg) 
= \int_{\domain[i]} \specificVolume[i] \source[\pressure] \dx,
\end{aligned}
\end{gather}
with the third term accounting for changes in fracture volume.

The dimension reduction is now performed for the energy balance, which for an equi-dimensional domain reads
\begin{gather}\label{eq:energy_balance}
\begin{aligned}
\dD{}{\timet}\int_{\domainEq}\density[f] \heatCapacity[f] \temperature \dx
 +\int_{\boundaryEq} \left(\advectiveHeatFlux + \conductiveHeatFlux \right) \cdot \dA = \int_{\domainEq}  \source[\temperature]\dx,
\end{aligned}
\end{gather}
where $\advectiveHeatFlux$ and $\conductiveHeatFlux$ are given by Eqs.~\eqref{eq:advectiveflux} and \eqref{eq:conductiveflux}. 
Utilising Eq.~\eqref{eq:dimension_reduction_flux}, the dimension reduction of the flux terms is
\begin{gather}\label{eq:heat_flux_tangential_normal}
\begin{aligned}
 \int_{\boundaryEq}\bigg(\advectiveHeatFlux + \conductiveHeatFlux\bigg)\cdot\dA 
 =\int_{\boundary{i}{}}\specificVolume[i]\bigg(\advectiveHeatFlux[i] + \conductiveHeatFlux[i]\bigg)\cdot\dA 
- \sum_{j\in \higherSet[i]} \projectFromInterface{l}{j}\bigg(\projectToInterface{h}{j}\specificVolume[h]\int_{\interface[j]} \interfaceAdvectiveHeatFlux[j] + \interfaceConductiveHeatFlux[j] \dx \bigg),
 \end{aligned}
 \end{gather}
and the internal boundary conditions on \boundary{h}{j} are 
\begin{gather}\label{eq:interface_Neumann_condition_heat}
    \begin{aligned}
    \advectiveHeatFlux[h] \cdot \normal[h]&=\projectFromInterface{h}{j}\interfaceAdvectiveHeatFlux[j], \\
    \conductiveHeatFlux[h] \cdot \normal[h]&=\projectFromInterface{h}{j}\interfaceConductiveHeatFlux[j].
	\end{aligned}
  \end{gather}
  The advective interface flux is defined according to the upstream direction of the interface fluid flux:
  \begin{gather}\label{eq:interface_advective_heat_flux}
    \begin{aligned}
            \interfaceAdvectiveHeatFlux = \left\{ \begin{array}{ l l }
      \interfaceFluidFlux \projectToInterface{h}{j} \density[f,h]\heatCapacity[f,h]\temperature[h] & \text{ if } \interfaceFluidFlux>0  \\
\interfaceFluidFlux \projectToInterface{l}{j} \density[f,l]\heatCapacity[f,l]\temperature[l] & \text{ if } \interfaceFluidFlux\leq0.
  \end{array} \right. 
    \end{aligned}
  \end{gather}
The Fourier-type conductive interface flux is
  \begin{gather}\label{eq:interface_conductive_heat_flux}
    \begin{aligned}
            \interfaceConductiveHeatFlux = - \heatConductivity[j]\frac{2}{\projectToInterface{l}{j}\aperture[l]} (\projectToInterface{l}{j} \temperature[l] - \projectToInterface{h}{j}\temperature[h]),
    \end{aligned}
  \end{gather}
  with the normal heat conductivity modelled as $\heatConductivity[j] = \projectToInterface{l}{j}\heatConductivity[f,l]$ since it originates from the dimension reduction of a fluid-filled domain. 
  
The dimension reduction of the remaining terms of Eq.~\eqref{eq:energy_balance} is performed by again invoking Eq.~\eqref{eq:dimension_reduction_scalar}, yielding
\begin{gather}\label{eq:energy_balance_linearised_simplified}
\begin{aligned}
\int_{\domain[i]} \heatCapacity[f]\density[f]\temperature[i]\dd{\specificVolume[i]}{\timet} 
&+ \specificVolume[i]\dd{}{\timet}\bigg(\density[f]\heatCapacity[f]\temperature[i]\bigg)  \dx 
+ \int_{\boundary{i}{}} \specificVolume[i]\bigg( \density[f]\heatCapacity[f]\temperature[i] \fluidFlux[i] 
- \heatConductivity[f] \nabla \temperature[i] \bigg)\cdot \dA \\
&- \sum_{j \in \higherSet[i]}\projectFromInterface{i}{j}\bigg(\projectToInterface{h}{j}\specificVolume[h]\int_{\interface[j]} \interfaceAdvectiveHeatFlux[j] + \interfaceConductiveHeatFlux[j] \dx \bigg) 
= \int_{\domain[i]} \specificVolume[i] \source[\temperature]\dx,
\end{aligned}
\end{gather}
where the second term in the first integral is decomposed similarly to Eq.~\eqref{eq:heat_balance_1term_2term}, with $\porosity=1$ in the calculation of effective coefficients.

\subsection{Fracture contact mechanics}\label{sec:fracture_deformation}
The traction balance, nonpenetration condition and friction law posed on a fracture $l$ are formulated in terms of interface displacements and fracture contact traction. The interface displacements on the two surfaces $\interface[j]$ and $\interface[k]$ are $\displacement[j]$ and  $\displacement[k]$, and the jump between the two sides is
\begin{gather}\label{eq:displacement_jump_definition}
  \begin{aligned}
\jump{ \displacement[l]} =  \projectFromInterface{l}{k}\displacement[k] - \projectFromInterface{l}{j}\displacement[j]. 
  \end{aligned}
\end{gather}
Since the fracture deformation depends on traction caused by the \textit{contact} between the two surfaces, the contribution from \pressure[l] should be subtracted on the fracture surfaces to yield the traction balance posed on the interfaces:
  \begin{gather}\label{eq:interface_traction_balance}
  \begin{aligned}
  \projectToInterface{l}{j}\traction[l]-\pressure[l] \identity\cdot \normal[l] &= \projectToInterface{h}{j} \stress[h] \cdot \normal[h]
  \quad \quad & \text{on }&\interface[j], \\
  \projectToInterface{l}{k}\traction[l]-\pressure[l] \identity\cdot \normal[l] &= - \projectToInterface{h}{k} \stress[h] \cdot \normal[h]
    \quad \quad & \text{on }&\interface[k].
    \end{aligned}
  \end{gather}
The fracture contact traction \traction[l] will for notational convenience be referred to as \traction in the following, and is defined according to the normal of the fracture, which is defined to equal \normal[h] on the $j$ side, i.e.\ $\normal[l]=\projectFromInterface{l}{j}\projectToInterface{h}{j}\normal[h]$. 
 When there is no mechanical contact between the interfaces, \traction is 0, implying that the higher-dimensional thermo-poromechanical tractions projected to the interfaces on the right-hand sides of Eq.~\eqref{eq:interface_traction_balance} are balanced by the fracture pressure.

A vector \is{\vect[]} defined on a fracture 
may be decomposed into the normal and tangential components
\begin{align}
\is{
    \iota_{n} = \vect[] \cdot \normal[l] \text{ and }   \vect[\tau] = \vect[] - \iota_{n} \normal[l].
}
\end{align}
With this notation, the nonpenetration condition reads
  \begin{gather}\label{eq:fracture_nonpenetration}
  \begin{aligned}
                \jump{\displacement}_n - \gap & \geq 0,  \\
                \normalTraction (\jump{\displacement}_n-\gap) & = 0,  \\
                \normalTraction & \leq 0,
\end{aligned}
\end{gather}
with the gap function \gap defined to equal the distance between the two fracture interfaces when in contact. The Coulomb friction law is
  \begin{gather}\label{eq:fracture_Coulomb}
  \begin{aligned}
                \norm{\traction[\tau]}& \leq -F\normalTraction,  \\
                 \norm{\traction[\tau]}& < -\frictionCoefficient \normalTraction \rightarrow \incrementTangentialDisplacement = 0, \\
                 \norm{\traction[\tau]}& = -\frictionCoefficient \normalTraction \rightarrow \exists \,\zeta \in \mathbb{R^+}:   \incrementTangentialDisplacement = \zeta\traction[\tau], 
\end{aligned}
\end{gather}
  with \frictionCoefficient denoting the friction coefficient and \incrementTangentialDisplacement denoting the tangential displacement increment.
In addition to enforcing the traction balance of Eq.~\eqref{eq:interface_traction_balance} and the conditions of Eqs. \eqref{eq:fracture_nonpenetration} and \eqref{eq:fracture_Coulomb}, a Dirichlet condition is assigned on \boundary{h}{j} so that
  \begin{equation}\label{eq:internal_boundary_displacement}
  \begin{aligned}
\projectFromInterface{h}{j}\displacement[j]=\displacement[h].
\end{aligned}
\end{equation}

\begin{figure}[tp]
\centering
\includegraphics[width=.95\textwidth]{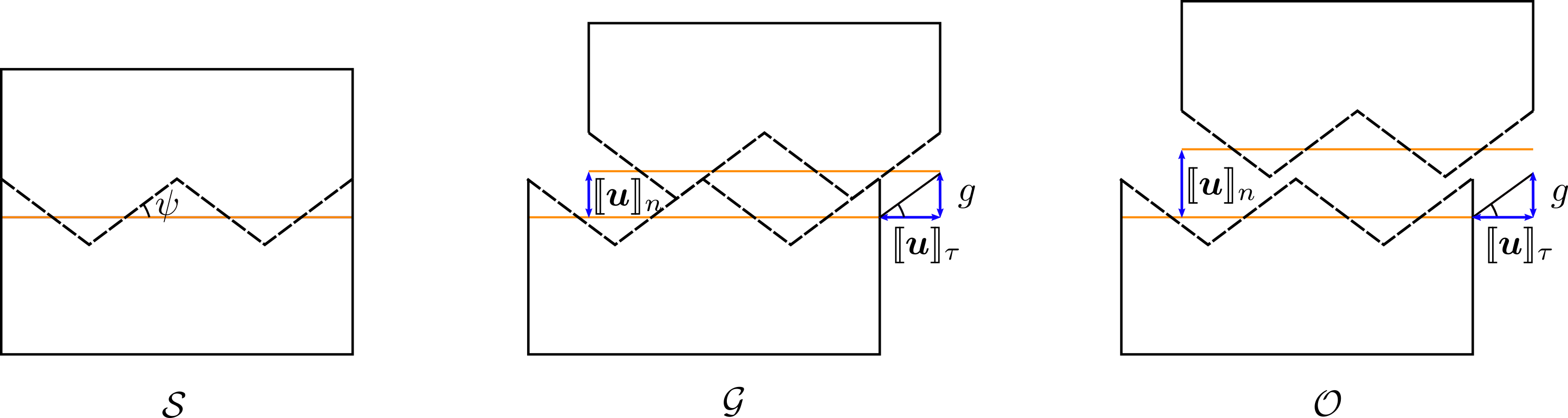}
\caption{Schematic representation of three fracture configurations: sticking (\stickingSet), gliding (\glidingSet) and open (\openSet). In the model, the fracture surfaces are represented as planar interfaces indicated by the orange lines. Idealised fracture roughness is shown by dashed sawtooth lines, with the inclination of the teeth equalling the dilation angle \dilationAngle, while the magnitude of displacement jumps and \gap are indicated by arrows.
In the first configuration, the fracture is undisplaced and closed with $g=\jump{\displacement[]}_n=\jump{\displacement[]}_\tangential=0$. In the second configuration, the fracture is still mechanically closed, but tangential displacement has resulted in fracture dilation due to roughness. In the third configuration, there is no mechanical contact across the fracture; that is, the fracture is mechanically open with $\jump{\displacement[]}_n > \gap$. 
}
\label{fig:dilation}
\end{figure}

The aperture introduced in Section \ref{sec:lower_dim_TH} is a function of displacement jump, $\aperture=\aperture\left(\jump{\displacement}\right)$. 
Due to roughness of the fracture surfaces, tangential displacements may induce dilation \cite{hossain2002shear} as illustrated in \figRef{dilation}.
The relationship between the dilation and the magnitude of tangential displacement is assumed to be linear and described by the dilation angle \dilationAngle following \citet{rahman2002shear_dilation}.
As modelled herein, the dilation is not merely a hydraulic effect impacting e.g.~the fracture permeability, but a mechanical effect in the sense that the normal distance between the fracture surfaces increases. As such, the dilation must be coupled back to the normal interface displacements and the matrix deformation through Eq.~\eqref{eq:internal_boundary_displacement}, which is achieved by choosing the gap function
  \begin{equation}\label{eq:gap_of_displacement}
  \begin{aligned}
  \gap= \tan(\dilationAngle)\norm{\jump{\displacement}_{\tangential}}.
  \end{aligned}
  \end{equation}
 The update is reversible; if the tangential displacement is reversed, \gap takes on its initial value. 
  
  Small-scale fracture roughness may provide a volume for the fluid to occupy even when the fractures are in an undeformed state. This leads to the following  relation between aperture and displacement:
\begin{equation}\label{eq:aperture_of_displacement}
\begin{aligned}
\aperture =  \aperture[0] + \jump{\displacement}_{n},
  \end{aligned}
  \end{equation}
where  \aperture[0] denotes the residual aperture in the undeformed state.

 In addition to entering the equations as a result of dimension reduction, \aperture governs the tangential permeability of a fracture or intersection line $i$ according to the cubic law \cite{zimmerman1996cubic},
 \begin{equation}\label{eq:cubic_law}
\begin{aligned}
\permeability[i]=  \frac{\aperture[i]^2}{12}\identity[i], 
  \end{aligned}
  \end{equation}
  where \identity[i] denotes the identity matrix of the fracture dimension. Equation \eqref{eq:cubic_law} constitutes a strongly nonlinear coupling, especially as \permeability[i] is multiplied by \specificVolume[] in Eq.~\eqref{eq:reduced_mass_balance_darcy}. 
Finally, the normal permeability of an interface is inherited from the lower-dimensional neighbour:
  \begin{equation}\label{eq:normal_permeability}
\begin{aligned}
\permeability[j]=\projectToInterface{l}{j} \permeability[l].
  \end{aligned}
  \end{equation}
  
\begin{figure}[tp]
\centering
\includegraphics[width=.6\textwidth]{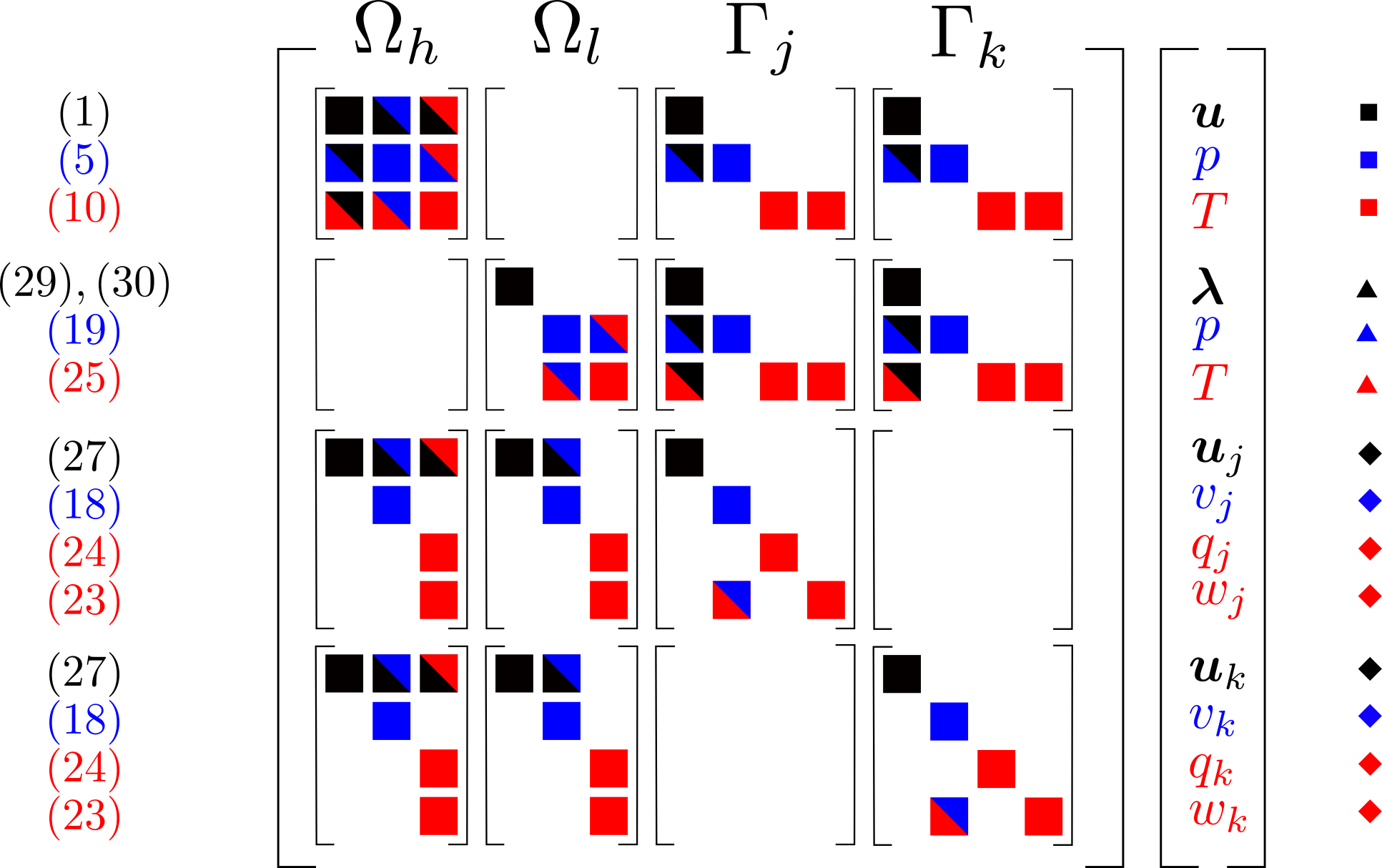}  \\
\hspace{10cm}
\includegraphics[width=.6\textwidth]{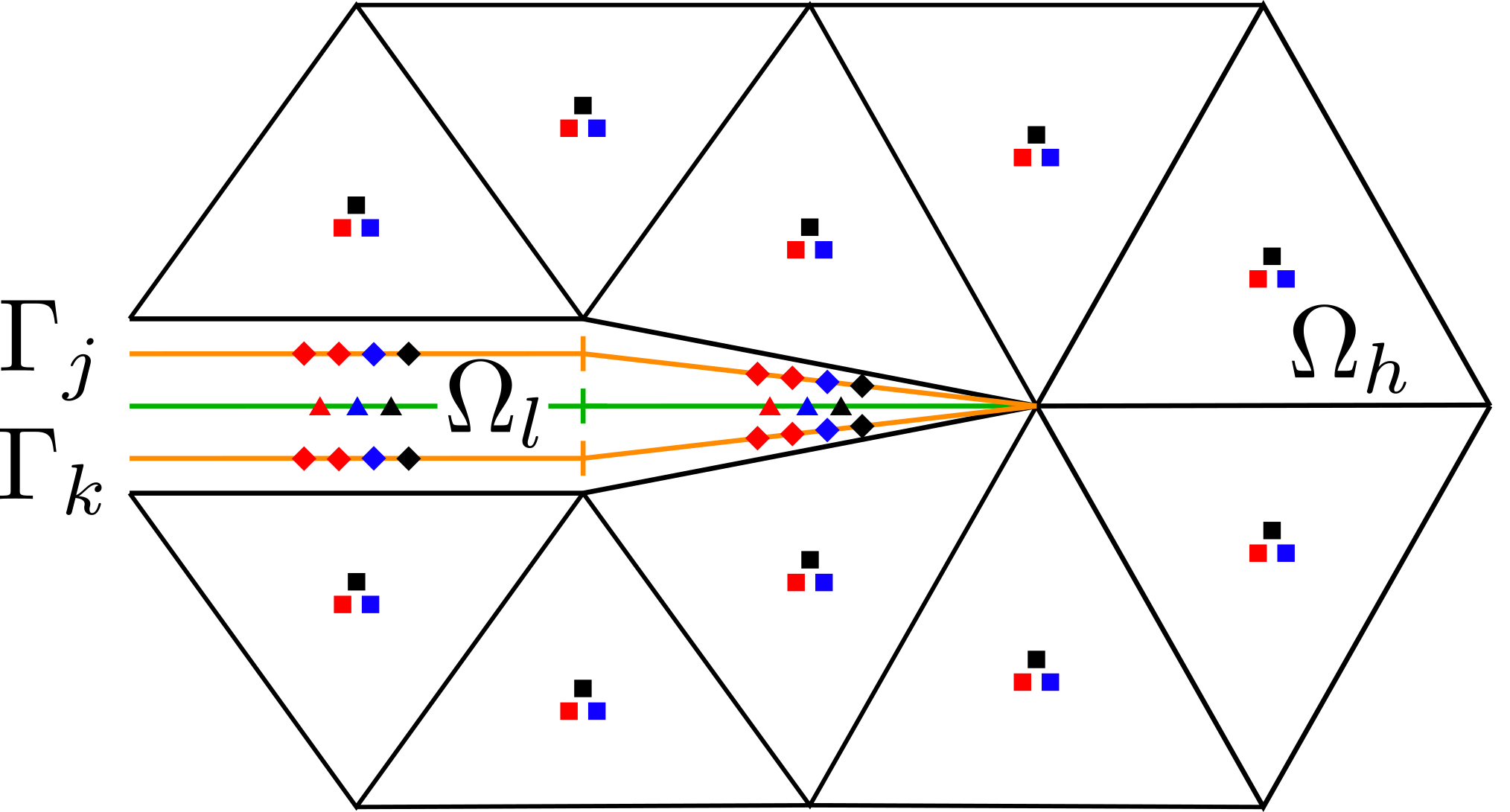}
\hfill
\includegraphics[width=.33\textwidth]{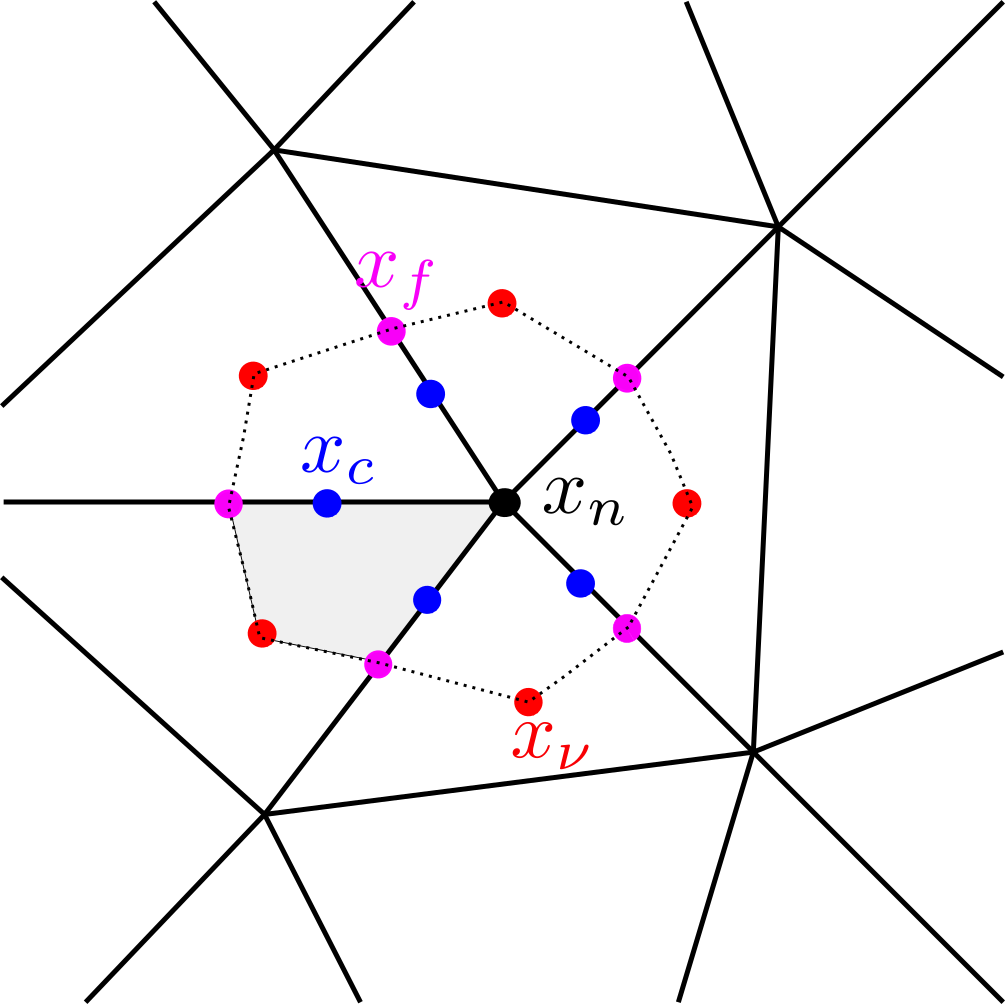}
\caption{Top: The two-level block structure represents the left-hand side of the linear equation system  for a matrix \domain[h], a fracture \domain[l] and interfaces \interface[j] and \interface[k] with corresponding equation numbers shown to the left. 
Bottom left: Spatial discretisation and spatial location of degrees of freedom for a domain corresponding to the equation system.  Matrix, fracture and interface grids are shown in black, green and orange, respectively, and the corresponding degrees of freedom are shown as squares, triangles and diamonds. Black represents displacement, blue pressure and red temperature; the relation between all markers and unknowns is shown at top right.
Bottom right: Subgrid around a node $x_n$ of the primary grid, which is shown in solid black lines. The O-shaped interaction region forming the stencil for the local systems is constructed by connecting the surrounding cell centres $x_\cell$ and face centres $x_f$ as indicated by the dotted lines. Continuity of primary variables is enforced in the points $x_c$; the shaded area indicates a subcell.
}
\label{fig:DOFs_blocks_and_mpxaregion}
\end{figure}

\section{Discretisation}\label{sec:discretisation}
This section describes the discretisation of the model presented in the previous section. 
The system is discretised in time using implicit Euler and solved monolithically using a direct solver \citep{umfpack}; a more scalable option would be to use iterative methods with block-preconditioners, following ideas in e.g.~\citep{franceschini2019block,both2019gradient}. 
The spatial grids are simplicial, and are
constructed such that the lower-dimensional cells coincide with
higher-dimensional faces; grids are generated by Gmsh
\cite{gmsh2009}.  
The model is implemented in the open-source fracture simulation toolbox PorePy presented \is{by} \citet{keilegavlen2020porepy}.

The mixed-dimensional framework gives rise to a two-level block
structure as the equations are discretised. 
The outer level corresponds to the subdomains and interfaces, with
entries internal to the subdomains on the diagonal and entries for the interdimensional
coupling on the off-diagonals. 
The inner level corresponds to the primary variables, 
with coupling effects between different variables on the
off-diagonals. The block structure is illustrated in \figRef{DOFs_blocks_and_mpxaregion}, which will be used in the following description of discretisation of individual terms by referring to the block in row $r$ and column $c$ as \block{r}{c} with $r$ and $c$ ranging from 1 to 14.

\subsection{Matrix thermo-poromechanics}
The spatial discretisation of the diffusive terms of the balance equations is achieved using a family of cell-centred finite volume schemes.
The approach is based on the multi-point flux approximation (MPFA) 
\citep{Aavatsmark2002} defined for diffusive scalar  problems and the multi-point stress approximation (MPSA) for vector problems \cite{nordbotten2014mpsa} and their combination for THM problems \cite{nordbotten2016biot, nordbotten2020mpxa}.
The scheme is formulated in terms of discrete displacement (\nd vectors), pressure and temperature unknowns and is locally momentum, mass and \is{energy} conservative.

The scheme's construction is based on a subdivision of the spatial grid as illustrated in \figRef{DOFs_blocks_and_mpxaregion}, with the gradients of displacement, pressure and temperature defined as piecewise constants on the subdivision.
The fluxes of the conserved quantities momentum, mass and energy are discretised via Hooke's, Darcy's and Fourier's law,  respectively.
Continuity is enforced for traction and mass and energy fluxes over faces of the subgrid and for the primary variables in the continuity points $x_c$, leading to one local system for the node of the primary grid. Each local system is partially inverted to express gradients in terms of the cell-centre values in nearby cells. A global system is constructed by collecting for each cell all face fluxes as expressed in terms of the cell-centred primary variables. For details, see \citet{nordbotten2020mpxa}. 

The coupling between the three equations is achieved by using the thermo-poroelastic stress for the local traction balances, which directly yields the contributions \block{1}{2} and \block{1}{3} representing the scalar variables' effect on the momentum balance.
\block{2}{1} and \block{3}{1}, which represent the displacement effects on the scalar balances, are constructed by assembly of the discrete divergence based on the local systems for the displacement gradients.

The standard finite volume  implicit Euler discretisation is applied to all time derivatives; that is, both the TH coupling blocks  \block{2}{3} and \block{3}{2} and the accumulation terms of \block{2}{2} and \block{3}{3}.
The advective term of \is{Eq.}~\eqref{eq:heat_balance2} is discretised using a first-order
upwind scheme, i.e.\ the temperature flux between cells $k$ and $l$ is
\begin{gather}\label{eq:upwind_discretisation}
    \begin{aligned}
    \left(\density[f] \heatCapacity[f] \temperature \fluidFlux\right)_{k,l} = \left\{ \begin{array}{ l l }
    \heatCapacity[f]\fluidFlux[k,l] \temperature[k]\density[f,k] & \text{ if } \fluidFlux[k,l]>0  \\
    \heatCapacity[f]\fluidFlux[k,l] \temperature[l]\density[f,l] & \text{ if } \fluidFlux[k,l]\leq0, 
  \end{array} \right.
    \end{aligned}
  \end{gather}
  with the fluid flux from cell $k$ to cell $l$ \fluidFlux[k,l] and \density[f] computed from the solution at the previous iteration.

\subsection{Mixed-dimensional flow and heat transfer}\label{sec:mixed_dim_discr}
All terms of the scalar  equations for the lower-dimensional subdomains are discretised using lower-dimensional versions of the corresponding \nd-dimensional discretisations. 
For the Darcy and Fourier fluxes, this implies that we use the MPFA scheme, while the advective fluxes are again treated by first-order upwinding.
The interdimensional coupling relations are discrete analogues to Eqs. \eqref{eq:interface_fluid_flux}, \eqref{eq:interface_conductive_heat_flux} and \eqref{eq:interface_advective_heat_flux}. Thus, they involve reconstruction of \pressure and \temperature on \boundary{h}{j}, which we base on discretisation matrices pertaining to the MPFA discretisations. For the matching grids used herein, the discrete projections are straightforward bijective mappings between faces of \domain[h] and cells of \interface[j] (\projectToInterface{h}{j} and \projectFromInterface{h}{j}) and between the cells of \domain[l] and \interface[j] (\projectToInterface{l}{j} and \projectFromInterface{l}{j}).

The nonlinearities arising through the products involving \aperture[] and \specificVolume[] are solved iteratively within the Newton scheme for fracture deformation described below. Specifically, the time derivatives are computed as additional  right hand side terms based on values from the previous iterate and time step.
However, the linear volume-change terms in the fractures are coupled fully implicitly to \displacement[j] so that the contribution for each fracture is the jump between the neighbouring higher-dimensional interfaces, as illustrated by the off-diagonal blocks \block{5}{7}, \block{5}{11}, \block{6}{7} and \block{6}{11}. 
Densities are computed from the solution at the previous iteration.
In some simulations involving strong advection and high temperature gradients, the density dependence in the gravity term of Darcy's law may lead to oscillatory fluxes between Newton iterations.
This may result in convergence problems related to the upstream discretisation of the advective term.
In these situations, convergence was achieved by damping the updates of the fluid flux of the advective term.

\subsection{Fracture contact mechanics}
Fracture deformation discretisation is based on the approach presented by \citet{hueber2008paper} and \citet{wohlmuth2011} with the frictional contact problem formulated as a variational inequality.
The formulation is expanded to account for the $\jump{\displacement}_\tangential$ dependency of \gap. Deformation constraints are reformulated as complementary functions $\complementaryFunction{}=\complementaryFunction{}(\primaryVector)$, with \primaryVector being the unknowns. 
The constraints are imposed by solving $\complementaryFunction{}=0$ through application of
the semismooth Newton method 
\begin{equation}\label{eq:newton}
\begin{aligned}
\jacobian{}(\primaryVector^k)(\increment\primaryVector^k)=-\complementaryFunction{}(\primaryVector^k), 
\end{aligned}
\end{equation}
where $\increment\primaryVector^k=\primaryVector^{k+1} - \primaryVector^{k}$. Correspondingly, the increment of a function $f$ between successive iterations $k$ and $k+1$ is $\increment f(\primaryVector^k) = f(\primaryVector^{k+1}) - f(\primaryVector^{k})$. $\jacobian{}$ is the generalised Jacobian of \complementaryFunction{}, i.e.\ the convex hull of the standard Jacobian wherever \complementaryFunction{} is differentiable. 

To facilitate imposition of different constraints based on the deformation states defined in Eqs.~\eqref{eq:fracture_nonpenetration} and \eqref{eq:fracture_Coulomb}, three disjoint sets describing the deformation state as open, sticking or gliding are defined:
\begin{gather}\label{eq:deformation_state_sets}
    \begin{aligned}
    \openSet&= \left\{\frictionBound \leq 0 \right\}, \\
    \stickingSet&= \left\{||\traction[\tangential]+\cNum\incrementTangentialDisplacement||<\frictionBound \right\}, \\
    \glidingSet &=\left\{||\traction[\tangential]+\cNum\incrementTangentialDisplacement||\geq \frictionBound >0  \right\}.
    \end{aligned}
  \end{gather}
Here, \cNum denotes a numerical parameter, the friction bound is
$\frictionBound=-\frictionCoefficient(\normalTraction + \cNum\left(\jump{\displacement}_{n}
- \gap\right))$ and \incrementTangentialDisplacement denotes the increment from the previous time step. 
Replacing \incrementTangentialDisplacement  by $\jump{\displacement[]}_\tangential$ in the above definition yields the cumulative fracture state sets, which are denoted by subscript $c$.

  The normal and tangential complementary functions are
\begin{gather}\label{eq:normal_complementary_function}
    \begin{aligned}
    \complementaryFunction{n}\left(\jump{\displacement}_{n}, \normalTraction\right) = - \normalTraction - \frac{1}{\frictionCoefficient} \maximum{0}{\frictionBound}
    \end{aligned}
  \end{gather}
  and
\begin{gather}\label{eq:tangential_complementary_function}
    \begin{aligned}
    \complementaryFunction{\tangential}\left(\jump{\increment\displacement}_{\tangential},\jump{\displacement}_{\tangential}, \traction{\tangential}\right) = \maximum{\frictionBound}{||\traction[\tangential] + \cNum\jump{\displacement}_{\tangential}||}\left(-\traction[\tangential]\right) + \maximum{0}{\frictionBound}(\traction[\tangential] +\cNum\incrementTangentialDisplacement).
    \end{aligned}
  \end{gather}
The corresponding generalised Jacobians are
\begin{gather}\label{eq:normal_jacobian}
    \begin{aligned}
    \jacobian{n}\left(\jump{\displacement}, \normalTraction\right) \left(\increment\jump{\displacement}, \increment\normalTraction\right) 
    &= - \increment\normalTraction - \characteristicFunction{\closedSet} \frac{1}{\frictionCoefficient}\increment\frictionBound 
    \end{aligned}
  \end{gather}
  and
  \begin{gather}\label{eq:tangential_jacobian_Hueeber_style}
    \begin{aligned}
    \jacobian{\tangential}\left(\jump{\displacement}, \jump{\dot{\displacement}}_{\tangential}, \traction\right) \left(\increment\jump{\displacement}, \increment\jump{\dot{\displacement}}_{\tangential}, \increment\traction\right) 
    =&- \maximum{\frictionBound}{\norm{\traction[\tangential]+\cNum\jump{\dot{\displacement}}_{\tangential}}}\increment\traction[\tangential] \\
    &-\characteristicFunction{\openSet\cup\glidingSet} \frac{\traction[\tangential]\left(\traction[\tangential]+\cNum\jump{\dot{\displacement}}_{\tangential}\right)^T}{\norm{\traction[\tangential] + \cNum\jump{\dot{\displacement}}_{\tangential}}}\left(\increment\traction[\tangential] +\cNum\increment\jump{\dot{\displacement}}_{\tangential}\right)\\
    &+\characteristicFunction{\closedSet} \, \frictionBound\left(\increment\traction[\tangential]+\cNum\increment\jump{\dot{\displacement}}_{\tangential}\right)\\
    &-\characteristicFunction{\stickingSet} \,\increment\frictionBound\traction[\tangential]\\
    &+ \characteristicFunction{\closedSet} \,\increment\frictionBound(\traction[\tangential]+\cNum\jump{\dot{\displacement}}_{\tangential}).
    \end{aligned}
  \end{gather}
  Here, \characteristicFunction{\star} is the characteristic function of a set $\star$ for a fracture cell \cell,
  \begin{equation}\label{eq:characteristic_function}
\begin{aligned}
\characteristicFunction{\star} = 
\left\{ \begin{array}{ l l }
    1 & \text{ if } \cell \in \star  \\
    0 & \text{ if } \cell \notin \star,  \\
  \end{array} \right.
\end{aligned}
\end{equation}
 while the increment of the friction bound is
 
  \begin{equation}\label{eq:d_friction_bound}
\begin{aligned}
\increment\frictionBound = -\frictionCoefficient[\increment\normalTraction +  \cNum \left(\increment\jump{\displacement}_{n} - \dGap\increment \jump{\displacement}_{\tangential}  \right)],
\end{aligned}
\end{equation}
with \dGap denoting the derivative of \gap with respect to $\jump{\displacement}$.
Hence, sorting each cell according to Eq.~\eqref{eq:deformation_state_sets} and imposing Eq.~\eqref{eq:newton} results in the following constraints:
\begin{gather}\label{eq:deformation_constraints}
    \begin{aligned}
	   \traction^{\cell, k+1} &= \vectorFont{0} & \cell \in \openSet, \\	
	   \jump{\displacement^{\cell, k+1}}_n - \left(\dGap\right)^{\cell, k}\incrementTangentialDisplacement[\cell, k+1] &=\gap^{\cell, k} - \left(\dGap\right)^{\cell, k}\incrementTangentialDisplacement[\cell, k]& \cell \in \glidingSet \cup \stickingSet, \\
	   \jump{\dot{\displacement}^{\cell, k+1}}_{\tangential} - 
	   \frac{\frictionCoefficient\incrementTangentialDisplacement[\cell, k]}{\frictionBound^{\cell, k}} \normalTraction^{\cell, k+1} &= \incrementTangentialDisplacement[\cell, k]& \cell \in \stickingSet, \\   
	   \traction[\tangential]^{\cell, k+1} - L^{\cell, k}\jump{\dot{\displacement}^{\cell, k+1}}_{\tangential} + \frictionCoefficient \fractureV^{\cell, k} \normalTraction^{\cell, k+1} &= \fractureR^{\cell, k} + \frictionBound^{\cell, k}\fractureV^{\cell, k} & \cell \in \glidingSet. 
    \end{aligned}
  \end{gather}
  The coefficients \fractureL, \fractureV and \fractureR are functions of
$\jump{\displacement}_{\tangential}^{k}$ and $\traction^{k}$, and can
thus be computed from the previous iterate.
 For the exact expressions and
further details of the discretisation and implementation of the fracture
deformation equations, see \citet{berge2020hmdfm}.  

The effect of letting \gap depend on 
$\jump{\displacement}_\tangential$ only appears in the normal condition in the two terms involving
the derivative \dGap.
The two cases $\gap=0$ and
Eq.~\eqref{eq:gap_of_displacement} will be considered below. The former obviously gives $\dGap=0$ while the latter gives
\begin{equation}\label{eq:dgap}
\begin{aligned}
\dGap = 
\left\{ \begin{array}{ l l }
    \tan(\dilationAngle) \frac{\jump{\displacement}_\tangential^T}{\norm{\jump{\displacement}_\tangential}} & \text{ if } \quad \norm{\jump{\displacement}_\tangential}>0  \\
    0 \ & \text{ if } \quad \norm{\jump{\displacement}_\tangential} = 0, 
  \end{array} \right.
\end{aligned}
\end{equation}
which may be inserted into Eq.~\eqref{eq:deformation_constraints} to finally yield \block{4}{4}, \block{4}{7} and \block{4}{11} of \figRef{DOFs_blocks_and_mpxaregion}.

\section{Numerical results}\label{sec:results}
This section presents three sets of simulations aimed at demonstrating the model's representation of complex process--structure interactions. 
In the first example, a convergence study is presented and coupling mechanisms investigated.
The second example explores different modelling choices for the relationship between displacement jumps and apertures. 
Finally, the model is applied to a geothermal scenario with a pressure stimulation phase and long-term cooling during a production phase.
Run scripts for the example simulations and animations showing temporal evolution of the solutions may be found in a dedicated GitHub repository \cite{porepyRunScripts}.

\subsection{Example 1 --- Convergence study}\label{sec:ex1}
While different components of the implementation have been verified in previous studies \cite{keilegavlen2020porepy, berge2020hmdfm}, analytical solutions probing the full model presented herein are not available. The first example is designed as a validation of the implementation: 
Starting from a coarse grid of \num{398} 2d cells, \num{38} 1d cells and two 0d cells, a sequence of six grids is produced by nested conforming refinement.
The finest grid, which has \is{\num{407551}} 2d cells and a total of \num{1647254} unknowns, is used as the reference solution for the convergence study and forms the basis of the discussion of coupling mechanisms.

\begin{figure}
    \centering
    \includegraphics[width=.65\textwidth]{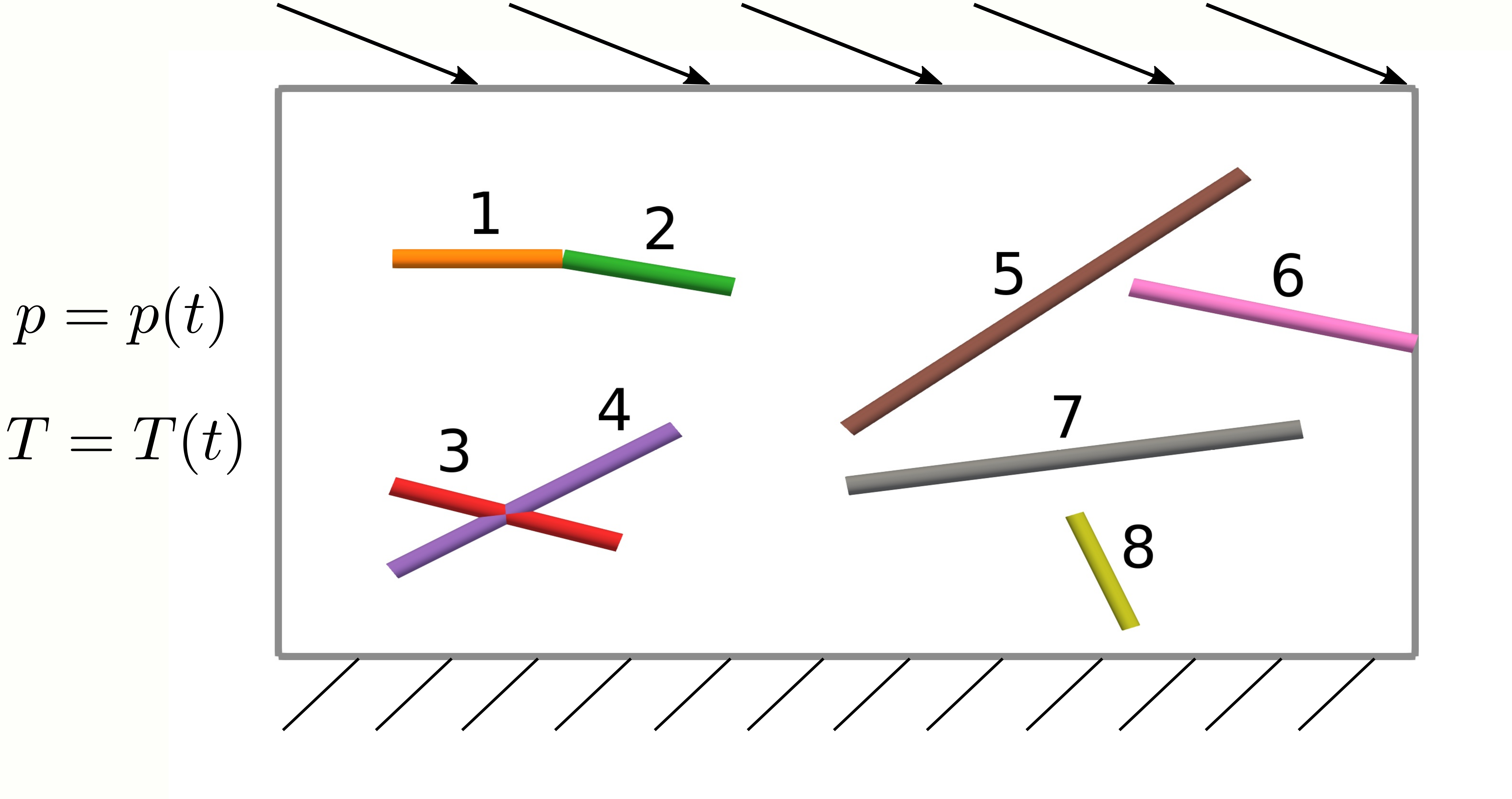}
    \caption{Fracture geometry and the boundary conditions driving the dynamics for examples 1 and 2. The colour scheme for the fractures is used throughout Sections \ref{sec:ex1} and \ref{sec:ex2}. The domain is fixed at the bottom and displaced at the top, while temporally varying pressure and temperature values $\pressure[\boundary{}{}]$ and \temperature[\boundary{}{}] are prescribed at the left boundary.} 
    \label{fig:ex1_geometry}
\end{figure}

The geometry of the two-dimensional domain with eight fractures is a modified version of a geometry presented in \citet{berge2020hmdfm} and is shown in \figRef{ex1_geometry}. 
It contains a kink formed by two fractures, an intersection formed by two other fractures and nearly intersecting fractures, as well as both immersed fractures and one fracture extending to the boundary.
These features can be expected to challenge the accuracy of numerical simulations.

Simulating three different phases allows us to distinguish between the influence of mechanical, hydraulic and thermal driving forces.
The three phases are defined through the boundary conditions as follows: Fixing the bottom and setting homogeneous stress conditions on the left and right boundary, a Dirichlet displacement value of $(\num{5e-4}, \num{-2e-4})^T$ \si{\metre} is applied at the top throughout the simulation and is the only driving force during phase I. Phase II begins when a pressure gradient of \SI{4e7}{\pascal} is applied from left to right. Once the solution has reached equilibrium, a boundary temperature \SI{15}{\kelvin} lower than the initial temperature is prescribed at the left boundary marking the onset of phase III. The initial values are \is{$\pressure=\SI{0}{\pascal}$, $\temperature =\temperature[0]=\SI{300}{\kelvin}$} and $\aperture[0]=\num{5e-4}$; no gravity effects are included in this example.

\begin{figure}
\centering
\includegraphics[width=1\textwidth]{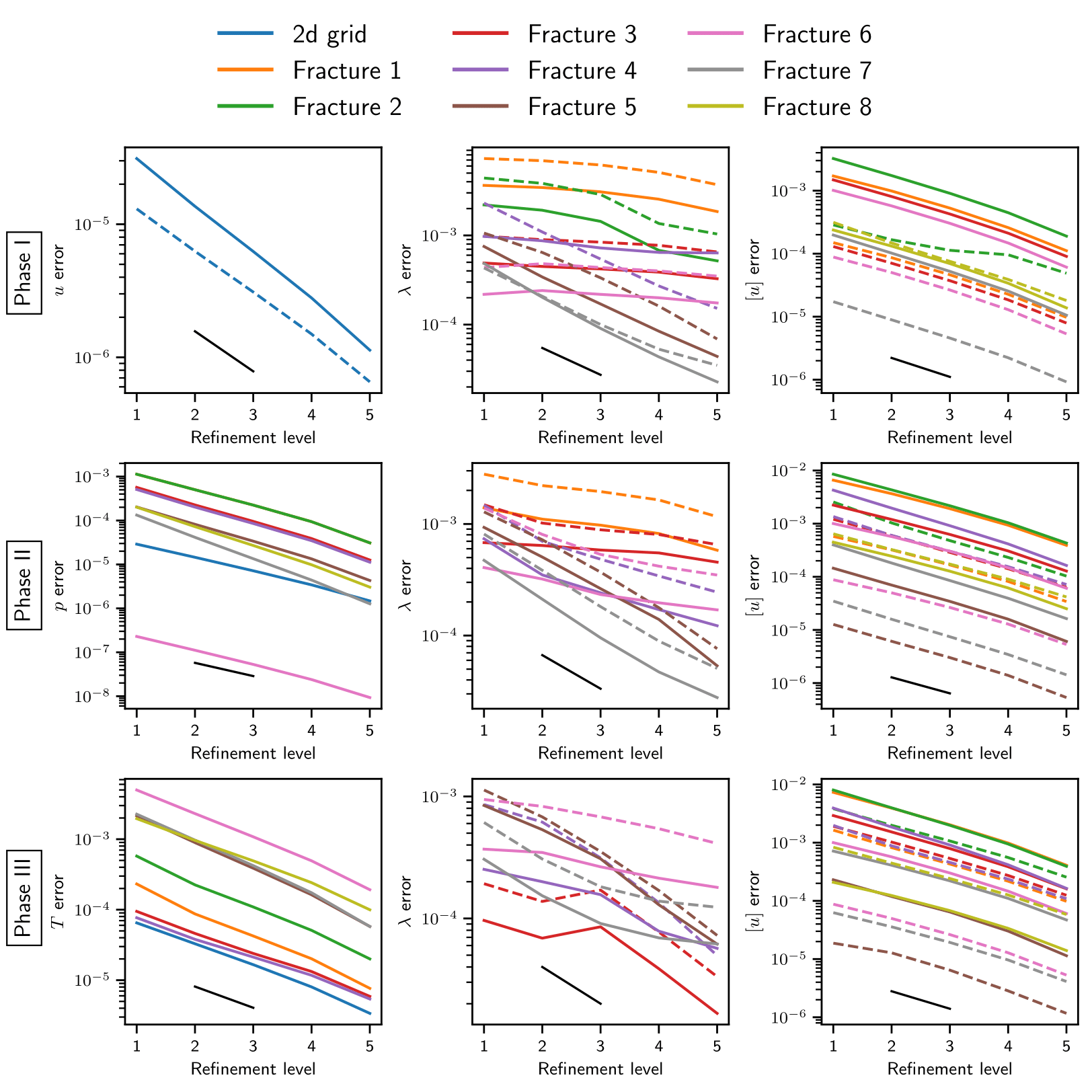}
\caption{Example 1: Errors relative to reference grid solution for solutions on five coarser grids at the end of the three phases, shown top to bottom. The three columns correspond to the variable of the main driving force, contact traction and displacement jumps. Solid and dashed lines correspond to $x$ and $y$ component in the matrix and tangential and normal component in the fractures. The black lines indicate first order.}
\label{fig:ex1_error_plots}
\end{figure}

\begin{figure}
\centering
\includegraphics[width=.48\textwidth]{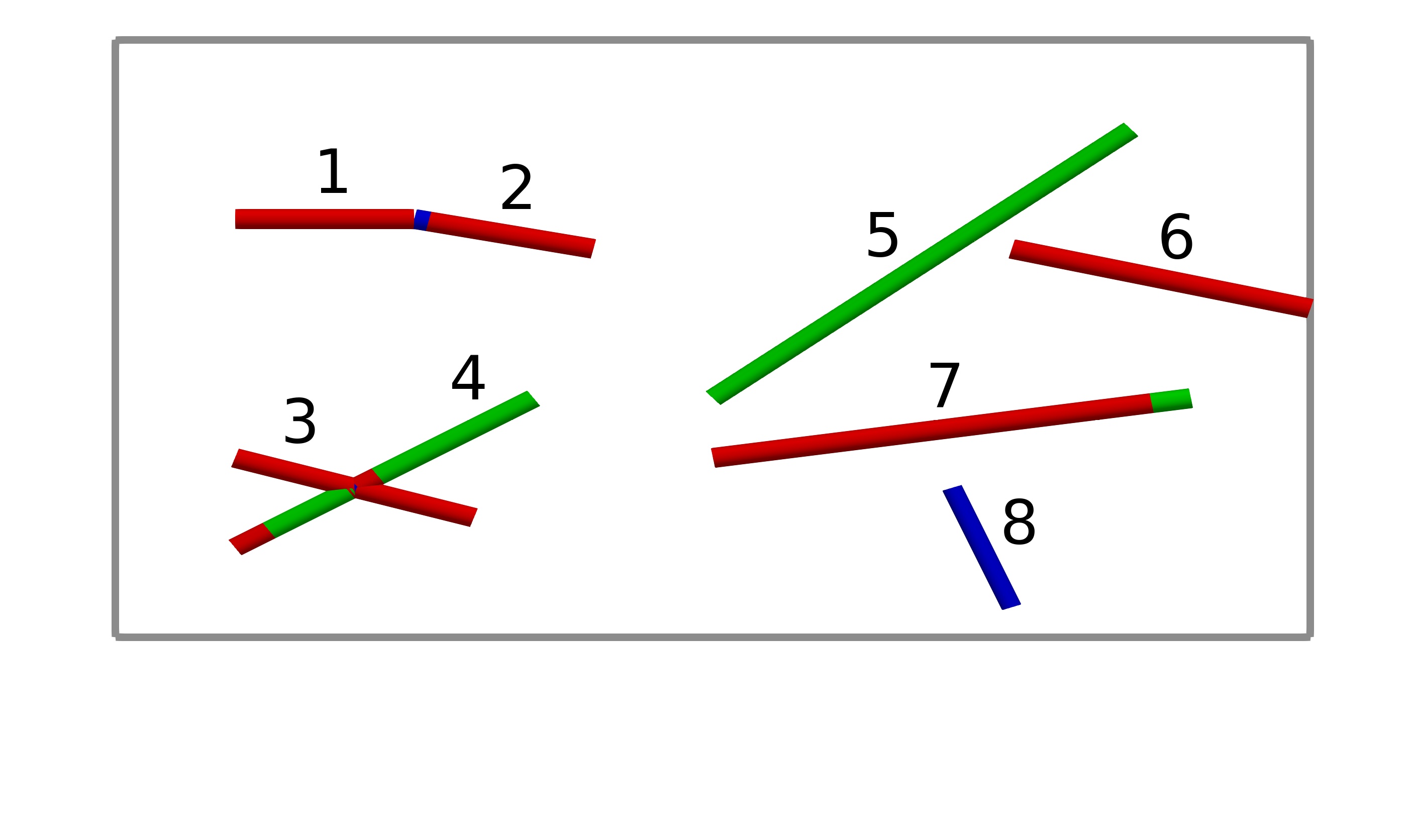}
\hfill
\includegraphics[width=.48\textwidth]{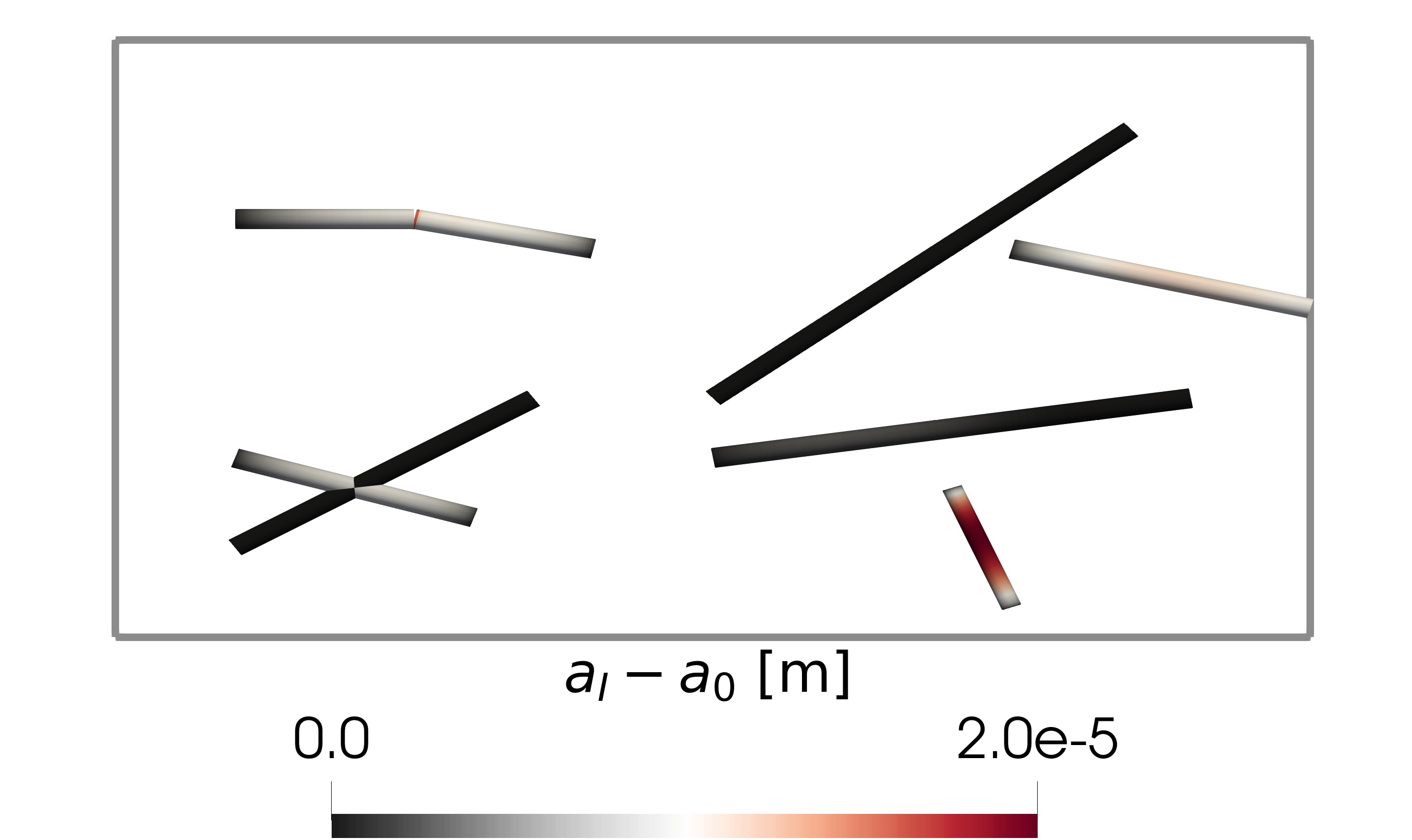}
\hfill
\includegraphics[width=.48\textwidth]{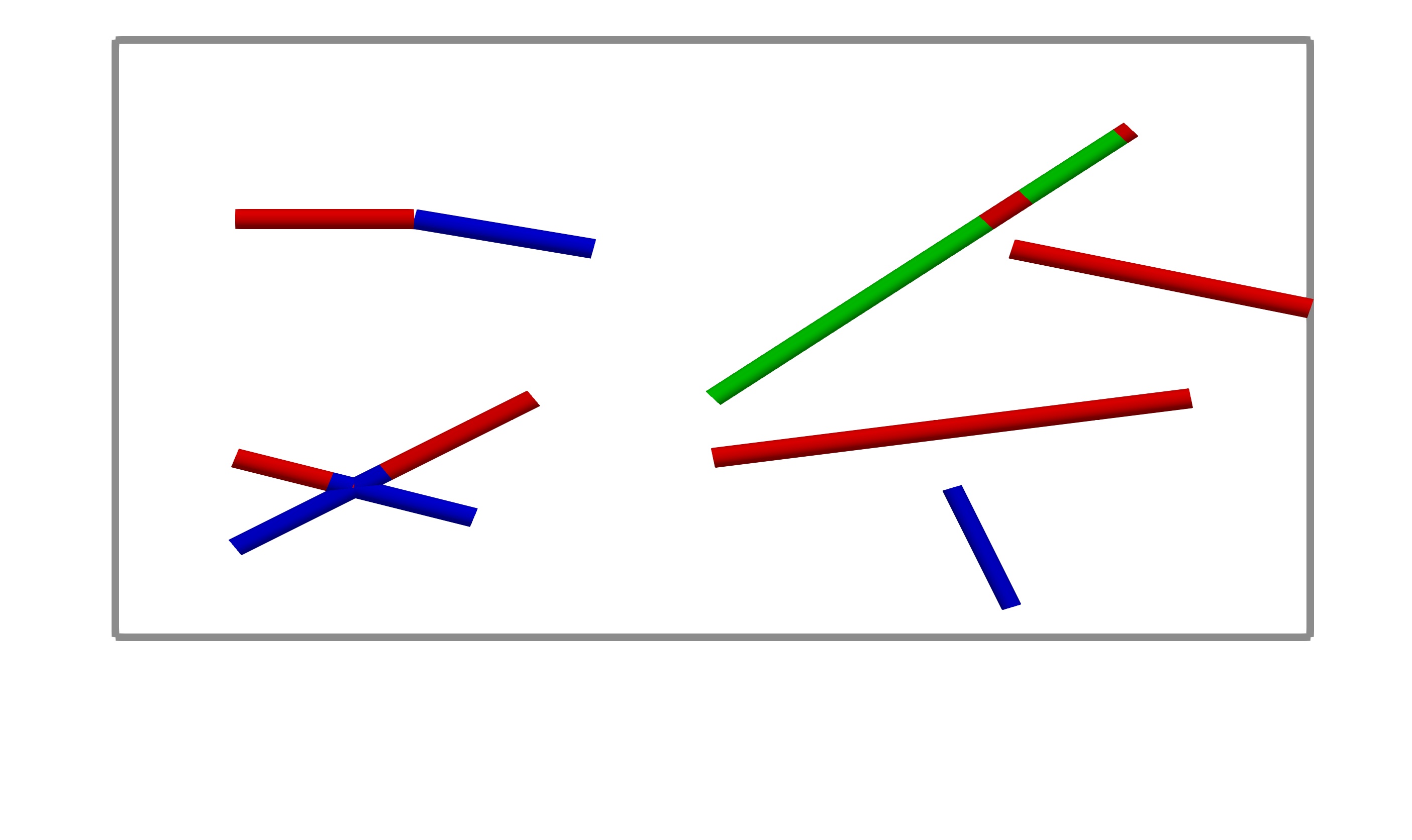}
\hfill
\includegraphics[width=.48\textwidth]{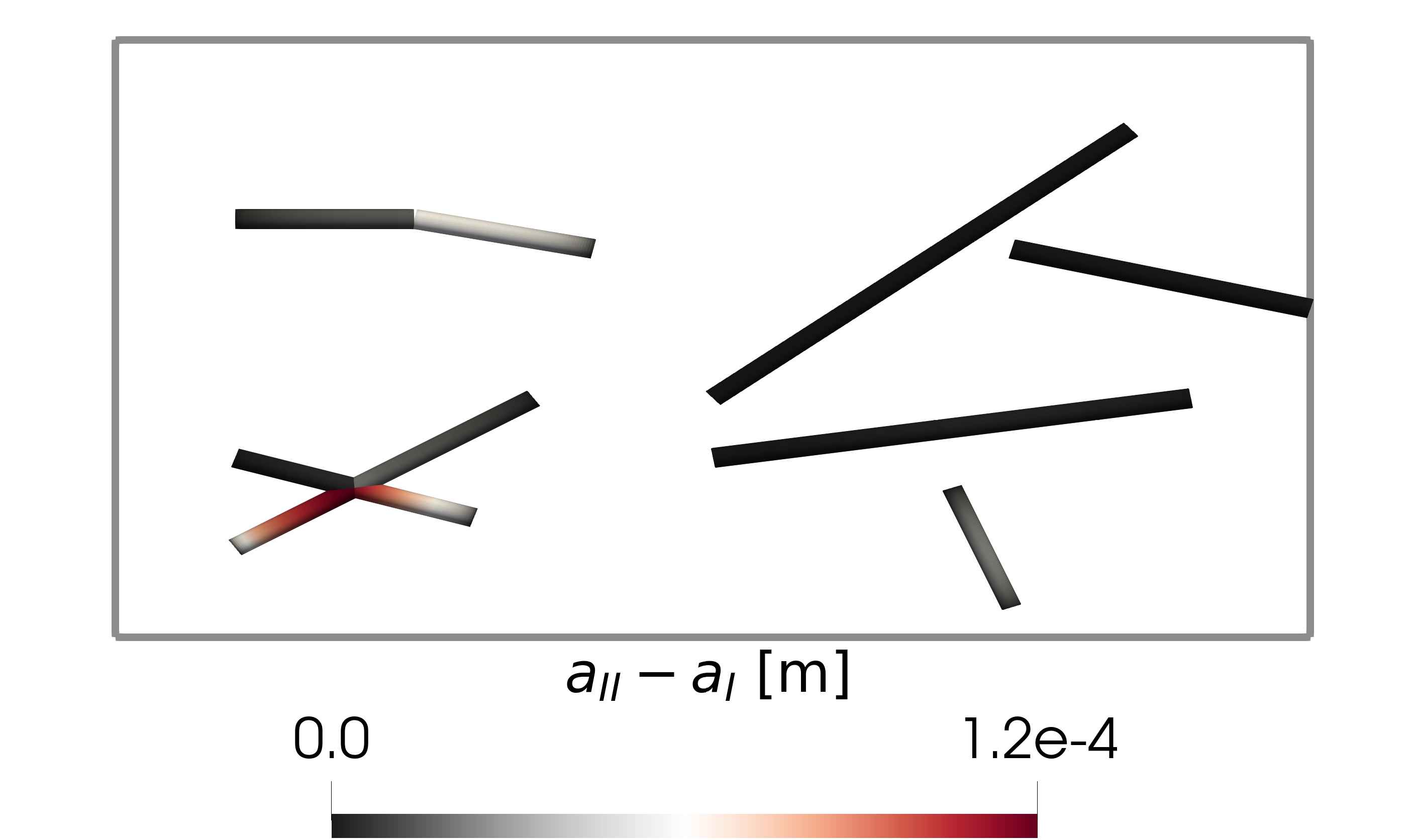}
\hfill
\includegraphics[width=.48\textwidth]{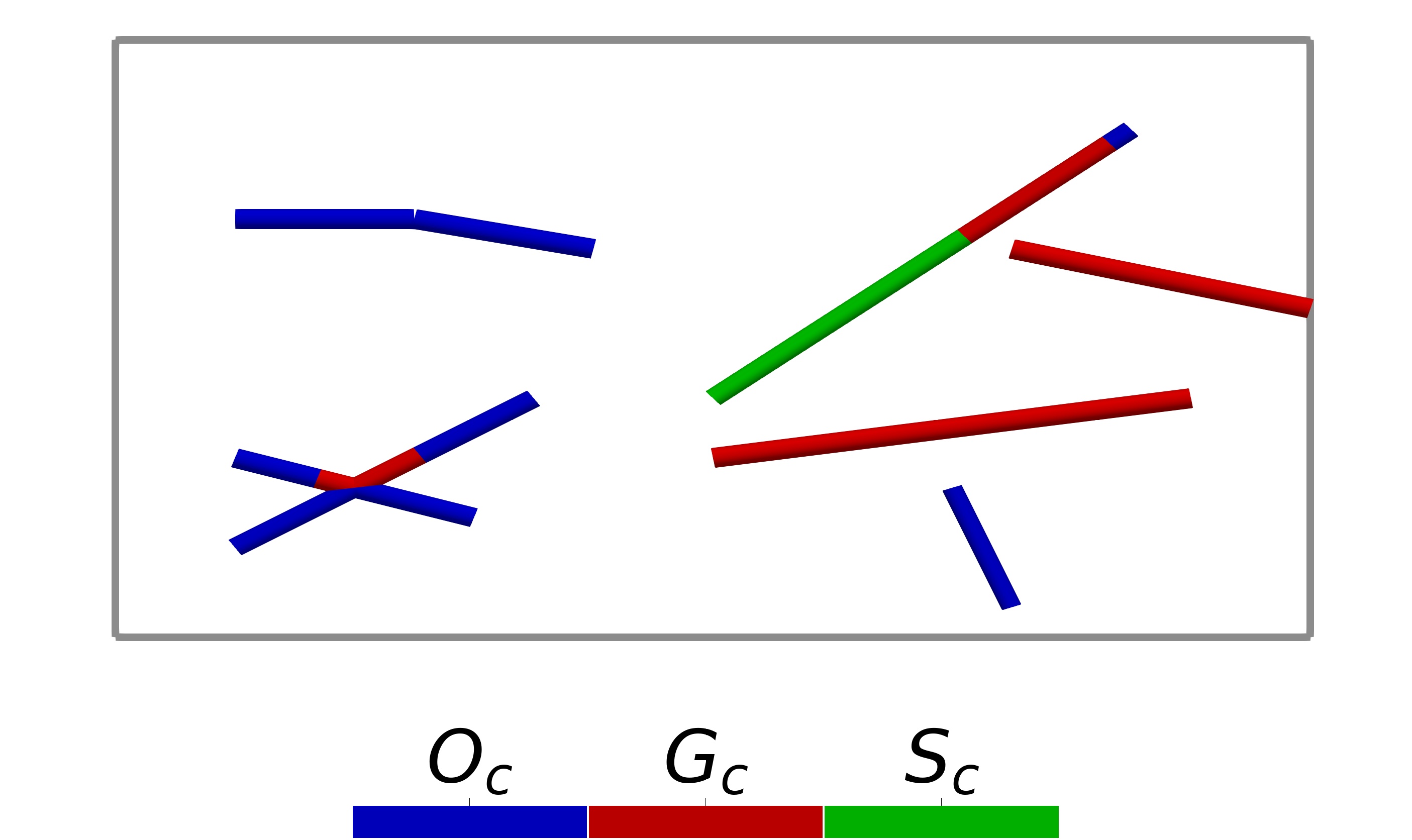}
\hfill
\includegraphics[width=.48\textwidth]{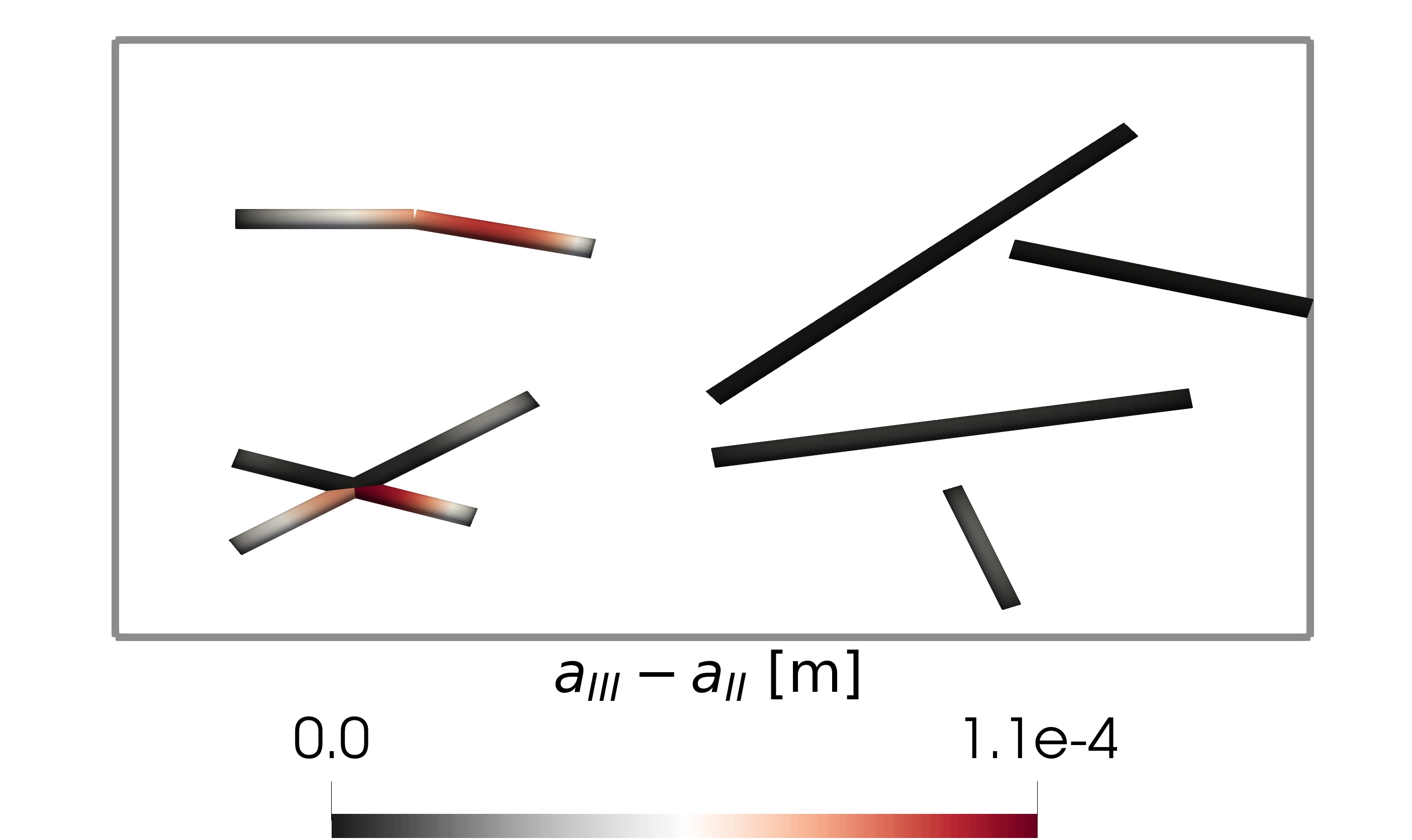} 
\caption{Example 1: 
Left: Fracture states according to cumulative displacement jumps at the end of phases I through III shown top to bottom. 
Right: Aperture increments throughout each of the three phases I through III shown top to bottom.
}
\label{fig:ex1_solutions}
\end{figure}
\begin{figure}
    \centering
    \includegraphics[width=.47\textwidth]{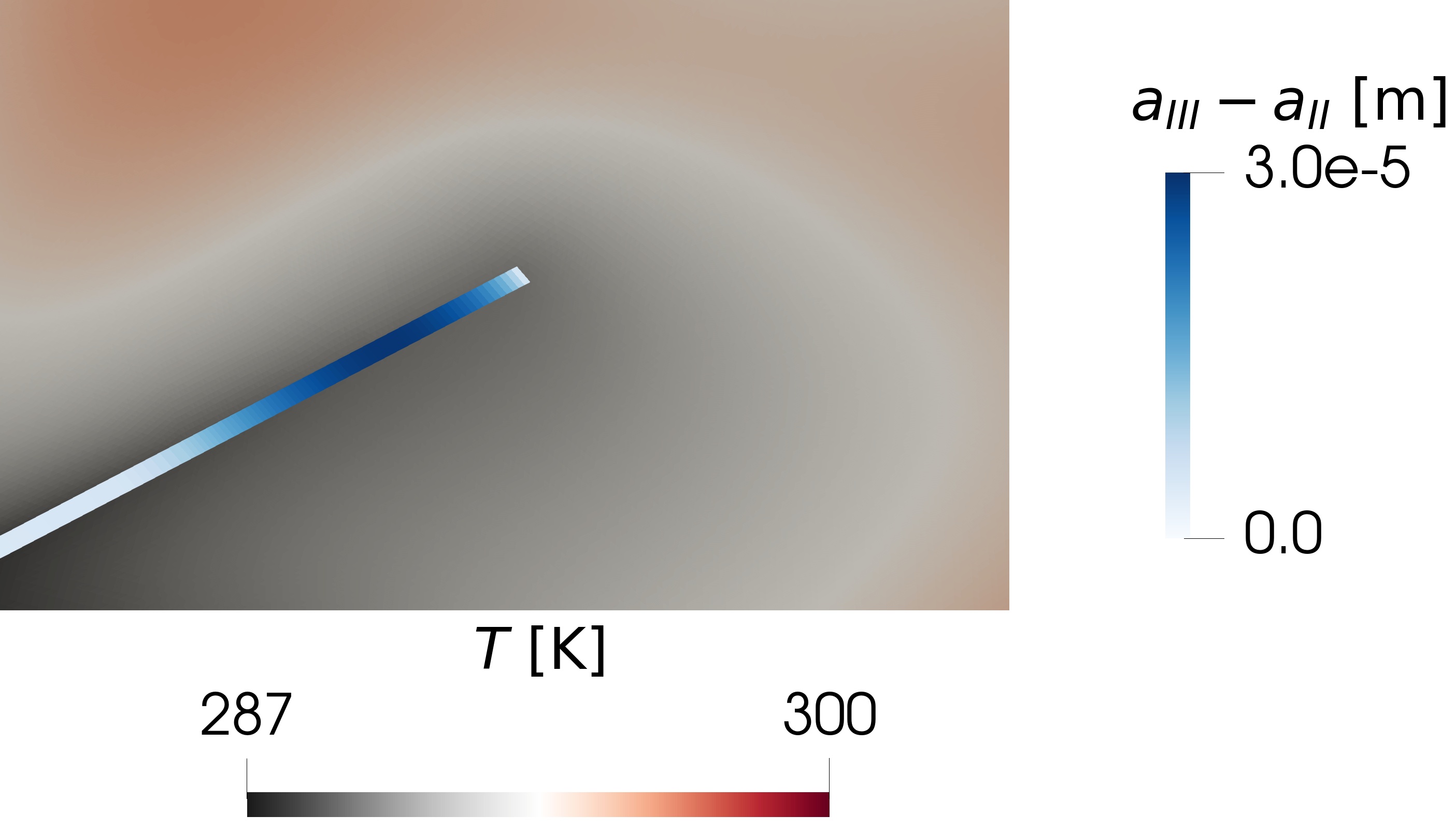}
    \hfill
    \includegraphics[width=.47\textwidth]{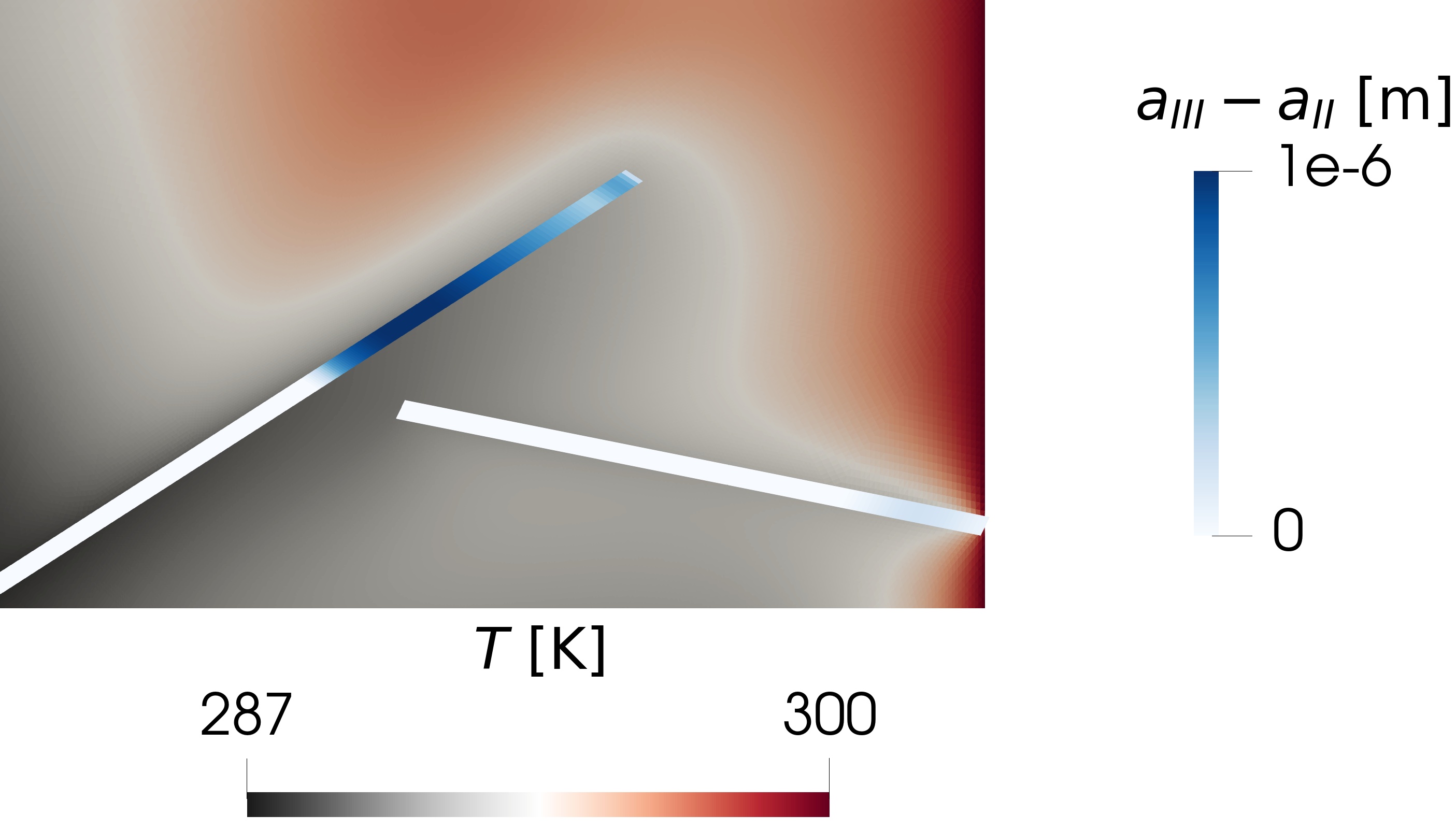}
    \caption{Example 1: 
    Matrix temperature with superimposed fracture aperture increment in the region surrounding fracture 4 (left) and fractures 5 and 6 (right).
    The temperature solutions are  \SI{8/3}{\hour} and \SI{22/3}{\hour} into phase III, at which point the cold temperature front has not yet moved past the respective regions.
The subscripts on \aperture indicate the value at the end of the corresponding phase, i.e.\ the values are the increments throughout phase III.}
    \label{fig:ex1_details}
\end{figure}

Figure \ref{fig:ex1_error_plots} shows convergence results for the end of the three phases. For each phase, we plot the errors for three primary variables on individual subdomains for the different refinement levels. The variables are displacement jumps, contact tractions and the variable related to the main driving force of the phase. The error is computed by projecting the cell-centre value of the coarse grids onto the reference grid and then computing the $L^2$ norm of the difference between coarse and fine solution. Errors are normalised by the number of reference cells in the subdomain multiplied by a weight $k$ representing the magnitude of the global range of the variable in question. The weights are obtained from the boundary conditions and are $k_\displacement=||(\num{5e-4}, \num{-2e-4})||$ \si{\metre}, $k_\pressure=$ \SI{4e7}{\pascal}, $k_\temperature = $\SI{15}{\kelvin} and $k_\traction = \text{E} k_\displacement$ with $\text{E}$ denoting Young's modulus.

In general, the expected first order convergence is observed. The exception is traction on some of the fractures (1, 2, 3 and 6). These local errors  may be attributed to the geometrical challenges posed by those fractures:
Fractures 1 and 2 meet in a kink, which seems to lead to relatively large errors compared to the remaining fractures as discussed by \citet{berge2020hmdfm}.
Fracture 3 intersects fracture 4, while the error on fracture 6 is concentrated around the leftmost tip, which is close to the neighbouring fracture 5.
However, the traction solutions converge for all fractures without small transition regions and the challenging geometrical features have no discernible effect on the convergence in the other primary variables.
 Therefore, taken together, the presented results serve as a verification of the model. 

Figure \ref{fig:ex1_solutions} shows the fracture deformation for each of the three phases, and thus demonstrates the effect of each of the three driving forces. The richness in physical processes and the complexity of coupling in the fractured THM problem is well illustrated by
a phenomenon observed towards the end of phase III around fractures 5 and 6 (see \figRef{ex1_details}). The role of fractures as preferential flow pathways leads to high flow rates and cooling in the region where fluid leaves fracture 5, both at the tip closest to the right boundary and in the area closest to fracture 6.
This, in turn, leads to local contraction of the matrix and fracture dilation --- in this particular case both through shear displacement and normal opening as seen from the final deformation state (bottom left in \figRef{ex1_solutions}).
The dilation further increasing the fracture conductivity can be expected to enhance the effect, which is also observed at the tip of fracture 4 somewhat earlier in the simulation.
This phenomenon of enhanced cooling-induced aperture increase in regions where the fluid enters or leaves a fracture can be expected to be of a general character. 
\begin{table}[ht]
  \centering
\begin{tabular}{|l  l  l |}
\hline
  \is{\bulkModulus{s}} & Solid bulk modulus & \SI{2.2e10}{\pascal} \\
  \is{\bulkModulus{f}} & \is{Fluid bulk modulus } & \is{\SI{2.5e9}{\per \pascal} }\\
  \shearModulus & Shear modulus & \SI{1.7e10}{\pascal} \\
  \viscosity & Viscosity & \SI{1.0e-3}{\pascal \second} \\
  \permeability & Permeability  & \SI{1.0e-15}{\metre \squared} \\
  \biotAlpha & Biot coefficient  & \num{0.8} \\\frictionCoefficient & Friction coefficient  & \num{0.5} \\
  \thermalExpansion[s] & Solid thermal expansion  & \SI{8.0e-6}{\per \kelvin} \\
  \thermalExpansion[f] & Fluid thermal expansion  & \SI{4.0e-4}{\per \kelvin} \\
  \heatConductivity[s] & Solid thermal conductivity  & \SI{3.0}{\watt \per \metre \per \kelvin} \\
  \heatConductivity[f] & Fluid thermal conductivity  & \SI{0.6}{\watt \per \metre \per \kelvin} \\
  \heatCapacity[s] & Solid specific heat capacity  & \SI{790}{\joule \per \kelvin} \\
  \heatCapacity[f] & Fluid specific heat capacity  & \SI{4.2e3}{\joule \per \kelvin} \\
  \porosity & Porosity  & \num{1.0e-2} \\
  \is{\density[s]} & Solid density & \SI{2.7e3}{\kilogram\per\metre} \\
  \density[f,0] & Reference fluid density & \SI{1.0e3}{\kilogram\per\metre} \\
\hline
\end{tabular}
\caption{Model parameters for the example simulations.}
\label{tab:parameters}
\end{table}

\subsection{Example 2 --- Fracture dilation models}\label{sec:ex2}
The second example is a study of different models for fracture dilation based on simulation of the case described in Section \ref{sec:ex1} with two simplified aperture models. 
In the first simplification, \noModel
, there is no coupling between shear displacement and dilation, i.e.\ $\gap=0$ and  $\aperture = \aperture[0] + \jump{\displacement}_n$. 
In the second simplified model (1) the aperture is related to the tangential displacement as $\aperture = \aperture[0] + \jump{\displacement}_n + \tan(\dilationAngle) \norm{\jump{\displacement}_\tangential} $ while \gap is kept constant.
This represents a naive one-way coupling which accounts for the dilation effect for the apertures and fracture permeability. 
We emphasise that neglecting the back-coupling to normal displacement -- and thus to the matrix momentum balance -- makes this model inconsistent. 
The model of \ref{sec:ex1},  where dilation is coupled to the displacement solution through the gap function according to Eq.~\eqref{eq:gap_of_displacement}, represents the full two-way dilation coupling and will be referred to as \strongModel. Thus, model names correspond to the number of directions of couplings accounted for by the models.

\begin{figure}
\centering
\includegraphics[width=1\textwidth]{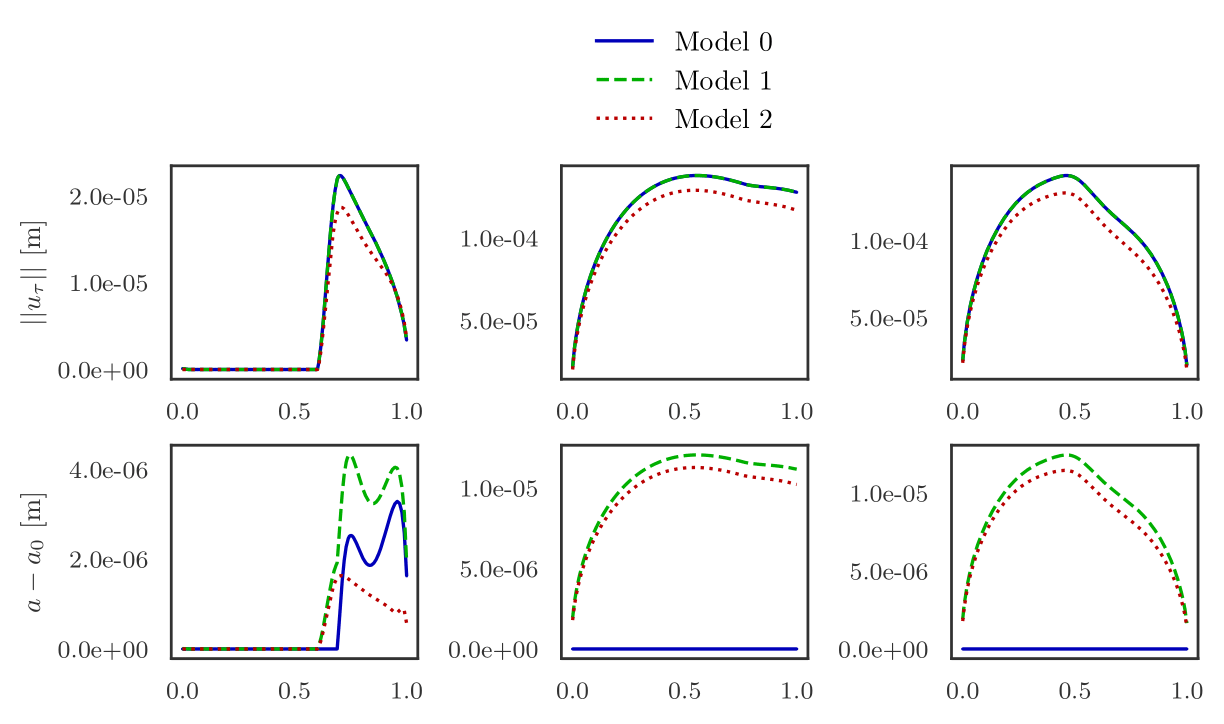}
\caption{Example 2: Final tangential displacement jumps (top) and apertures (bottom) along \is{the} three (partially) closed fractures \is{5, 6 and 7, corresponding to the three columns left to right. Results are shown for} the three different models for the relationship between \aperture and $\jump{\displacement}$. The cells are sorted from lowest to highest $x$ coordinate.}
\label{fig:fracture_aperture_and_tangential_jump_refinement_4}
\end{figure}

A comparison in terms of the final spatial distribution of aperture increase and tangential displacement jump on each of the closed fractures is shown in \figRef{fracture_aperture_and_tangential_jump_refinement_4}. 
As the dilation relations are irrelevant for open fractures, analysis is based on the mostly closed fractures  5 through 7.
Fractures 6 and 7 clearly demonstrate how the dilation coupling in \strongModel reduces tangential displacement compared to the simplified methods, as the induced normal displacement increases the normal traction on the fractures. 
Interestingly, the apertures displayed in \figRef{fracture_aperture_and_tangential_jump_refinement_4} show overestimation for the one-way coupling due to the above-mentioned overestimation of the tangential jumps.
Model 0 obviously yields no shear dilation. The \noModel aperture increase of fracture 5 is thus related to 
the fracture being open. Note that part of this region is closed for \strongModel (cf. the bottom left illustration of \figRef{ex1_solutions}) demonstrating how inconsistency affects the results beyond the prediction of \aperture.

The magnitude of the discrepancy must be expected to depend on rock and fault parameters, particularly the dilation angle.
The results demonstrate the strong coupling in the problem, and that compromising this coupling in numerical modelling may lead to significant error in the results.

\subsection{Example 3 --- Hydraulic stimulation and long-term cooling of a geothermal reservoir}\label{sec:ex3}
\begin{figure}[tp]
\centering
\includegraphics[width=.42\textwidth]{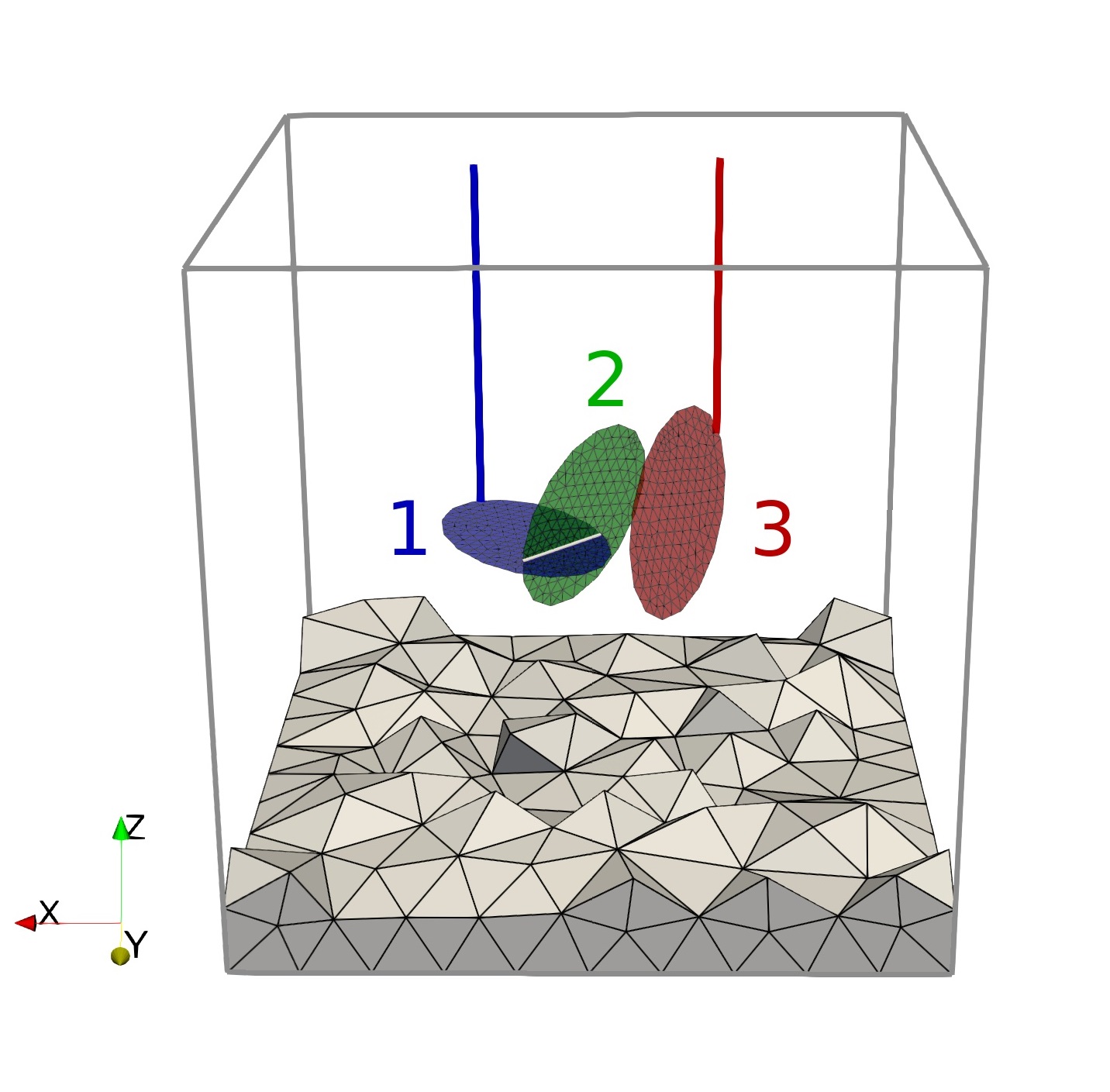}
\hfill
\includegraphics[width=.57\textwidth]{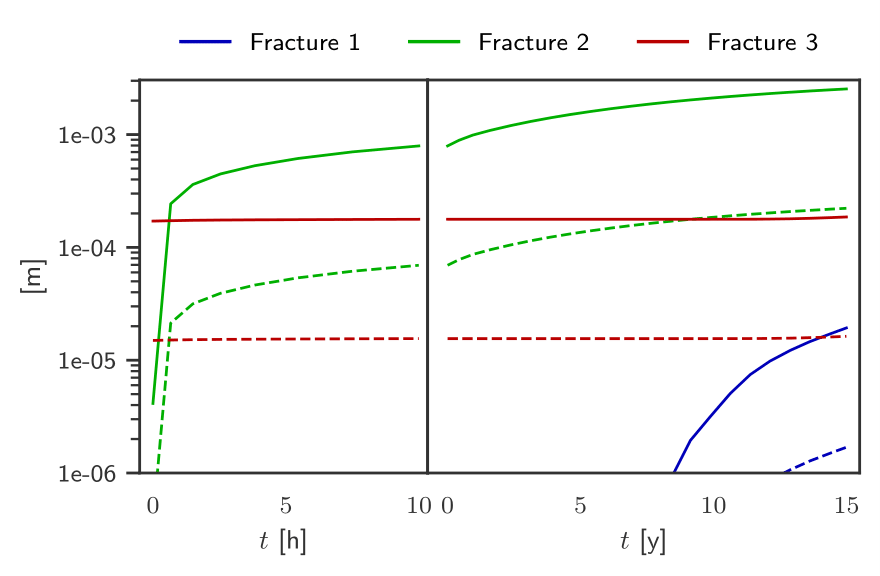}
\caption{Example 3: Left: Fracture network geometry and well paths. The grey lines indicate the domain boundary and the white line is the 1d fracture intersection,  while the injection and production wells are indicated by blue and red lines, respectively. Also shown are the 2d grid cells and some coarse 3d grid cells close to the boundary, indicating grid refinement in the region of interest.
Right: $L^2$ norm of tangential (solid lines) and normal (dashed lines) displacement jumps on each fracture during phases II and III. The values are normalised by the number of fracture cells. The black dashed line shows the number of Newton iterations needed for convergence. }
\label{fig:ex3_geometry_and_time}
\end{figure}

The third example shows hydraulic stimulation of a geothermal reservoir, followed by an injection and production phase leading to long-term reservoir cooling for the   three-dimensional geometry in \figRef{ex3_geometry_and_time}. The domain is the box $(\SI{-750}{\metre}, \SI{750}{\metre}) \times (\SI{-750}{\metre}, \SI{750}{\metre}) \times (\SI{-1750}{\metre}, \SI{-250}{\metre})$ and contains three fractures, two of which intersect along a line, and two wells.
The initial values are \is{$\pressure=\pressure[H]=\density[f,0]\gravity z$ \si{\pascal}, $\temperature=\temperature[0]=\SI{350}{\kelvin}$} and $\aperture[0]=\SI{2e-3}{\metre}$, with the positive direction of the $z$ axis pointing upwards and $\gravity=\SI{9.81}{\metre\per\second\squared}$ denoting the gravitational acceleration.
After letting the system reach equilibrium under the mechanical boundary conditions representing an anisotropic background stress in phase I, we simulate a pressure stimulation phase (II) and a production and long-term cooling phase (III). 
In the 10 hour stimulation phase, the flow rates of the injection and production wells are \num{75} and \SI{0}{\liter \per \second}, respectively.
During the 15 year production phase, both rates are \SI{20}{\liter \per \second}. The injection temperature is \SI{70}{\kelvin} below the reservoir temperature. 
The wells are incorporated as source terms in the fracture cells intersected by the well paths, with upwind discretisation for the energy equation in the production cell.
Hydrostatic Dirichlet boundary conditions $\pressure=\pressure[H]$ apply for the pressure. An anisotropic compressive background stress is imposed with the following non-zero stress tensor values
\begin{gather}\label{eq:ex3_stress_tensor}
\begin{aligned}
\stress[xx]=\frac{3}{4} \density[s] \gravity z \qquad \stress[yy]=\frac{3}{2} \density[s] \gravity z \qquad \stress[zz]=\density[s]\gravity z.
\end{aligned}
\end{gather}
While the remaining parameters listed in the \tabRef{parameters} are plausible for geothermal reservoirs, they do not correspond to a specific site. \is{The number of Newton iterations required for convergence is below ten for all time steps.}

The results are summarised through the temporal evolution of the norm of the displacement jumps on the three fractures shown in \figRef{ex3_geometry_and_time}. Significant stimulation effects appear in both \is{phase II and phase III}, with the magnitude of the jumps somewhat larger during the cooling; dynamics initiate latest on the injection fracture due to its orientation relative to the background stress.

\begin{figure}[tp]
\centering
\includegraphics[width=.42\textwidth]{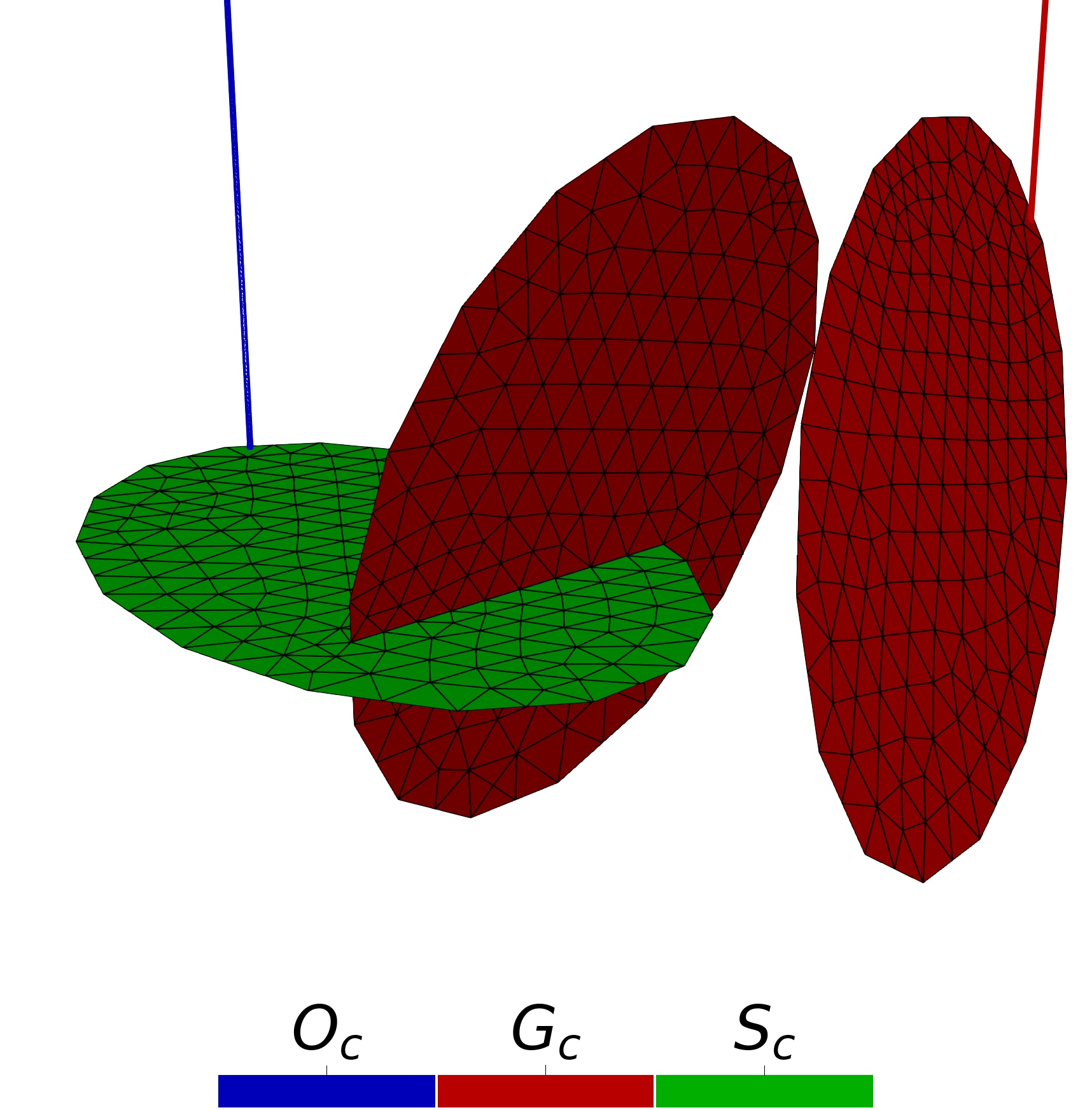}
\hfill
\includegraphics[width=.42\textwidth]{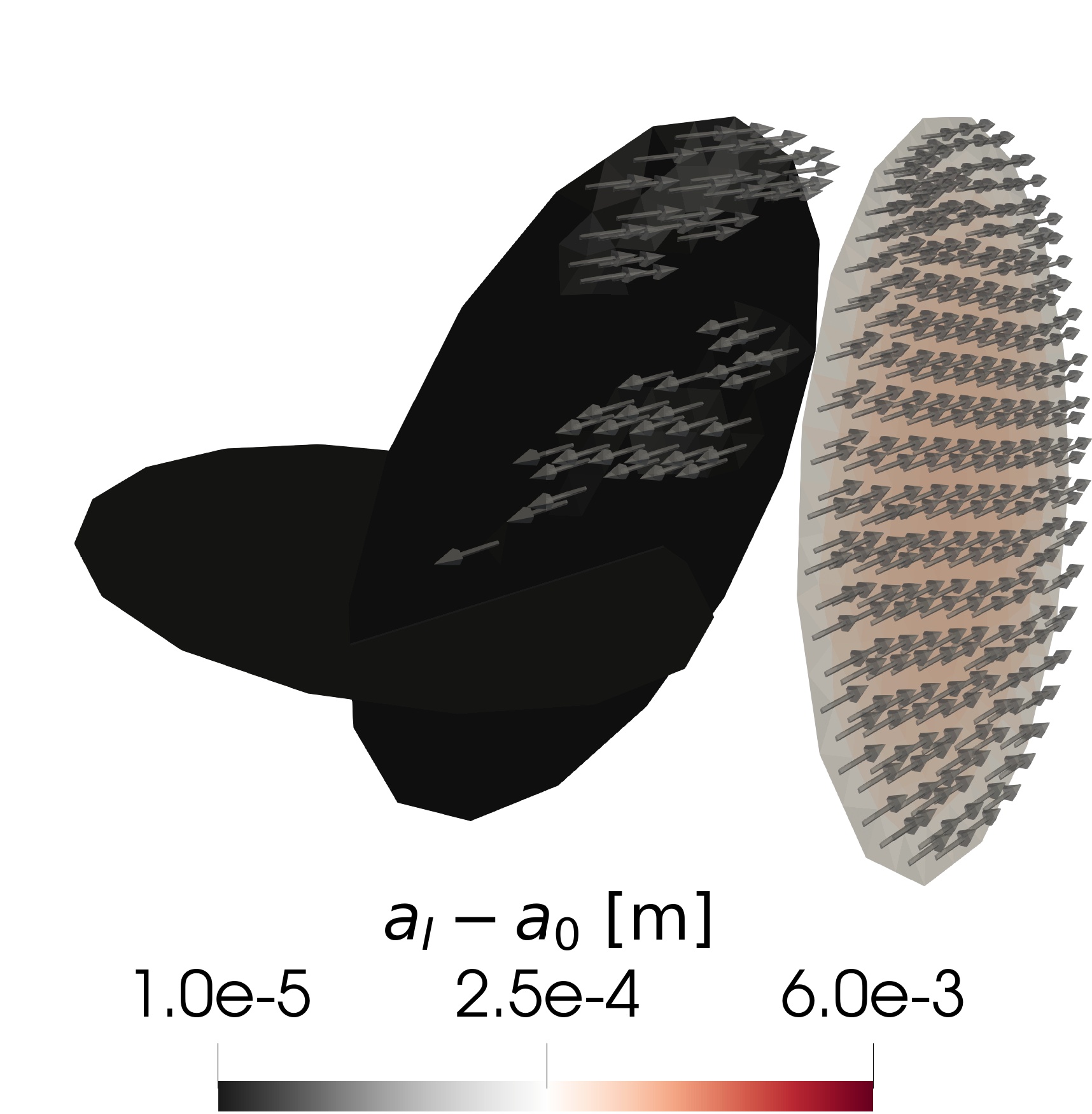}
\includegraphics[width=.42\textwidth]{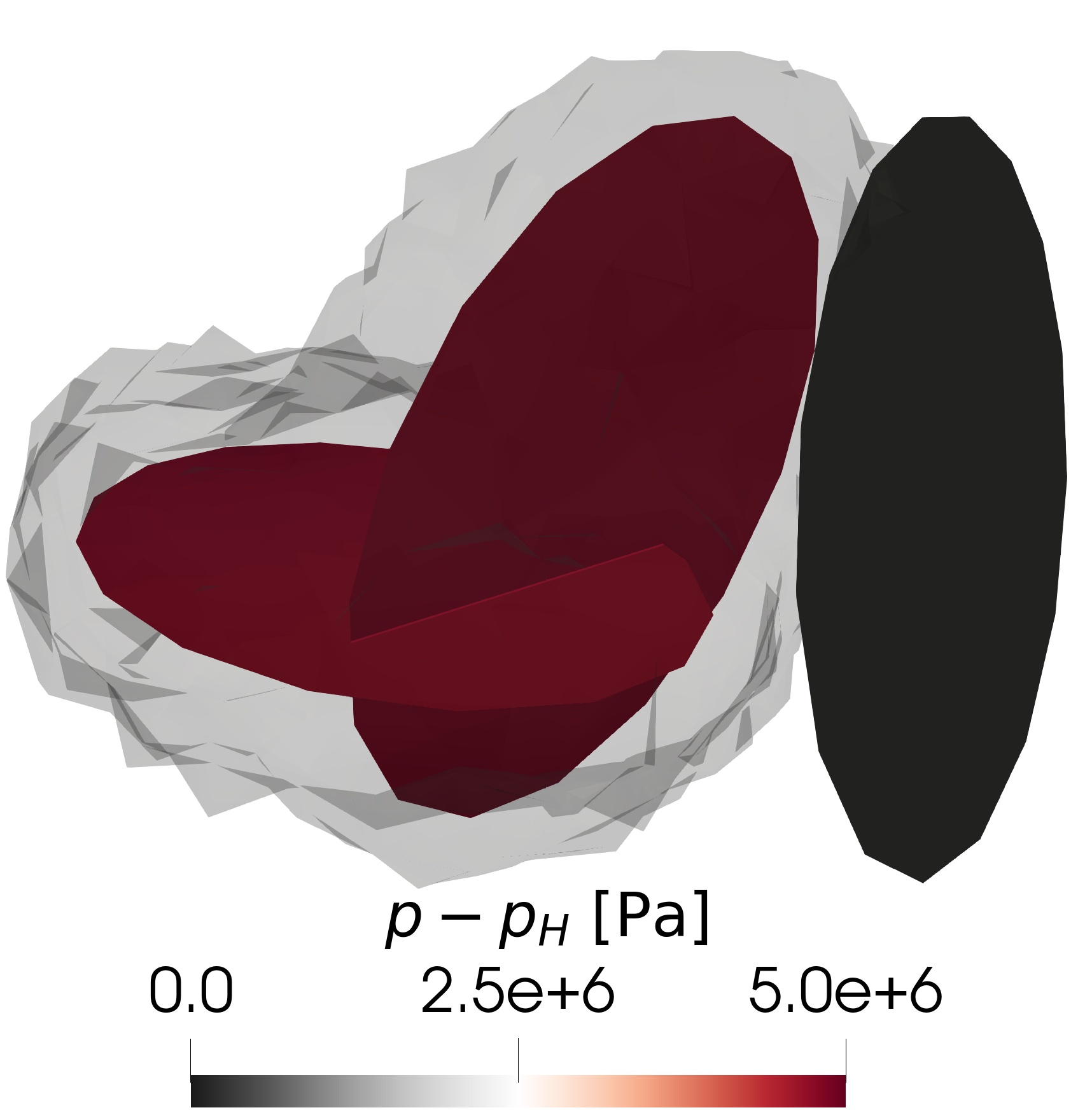}
\hfill
\includegraphics[width=.42\textwidth]{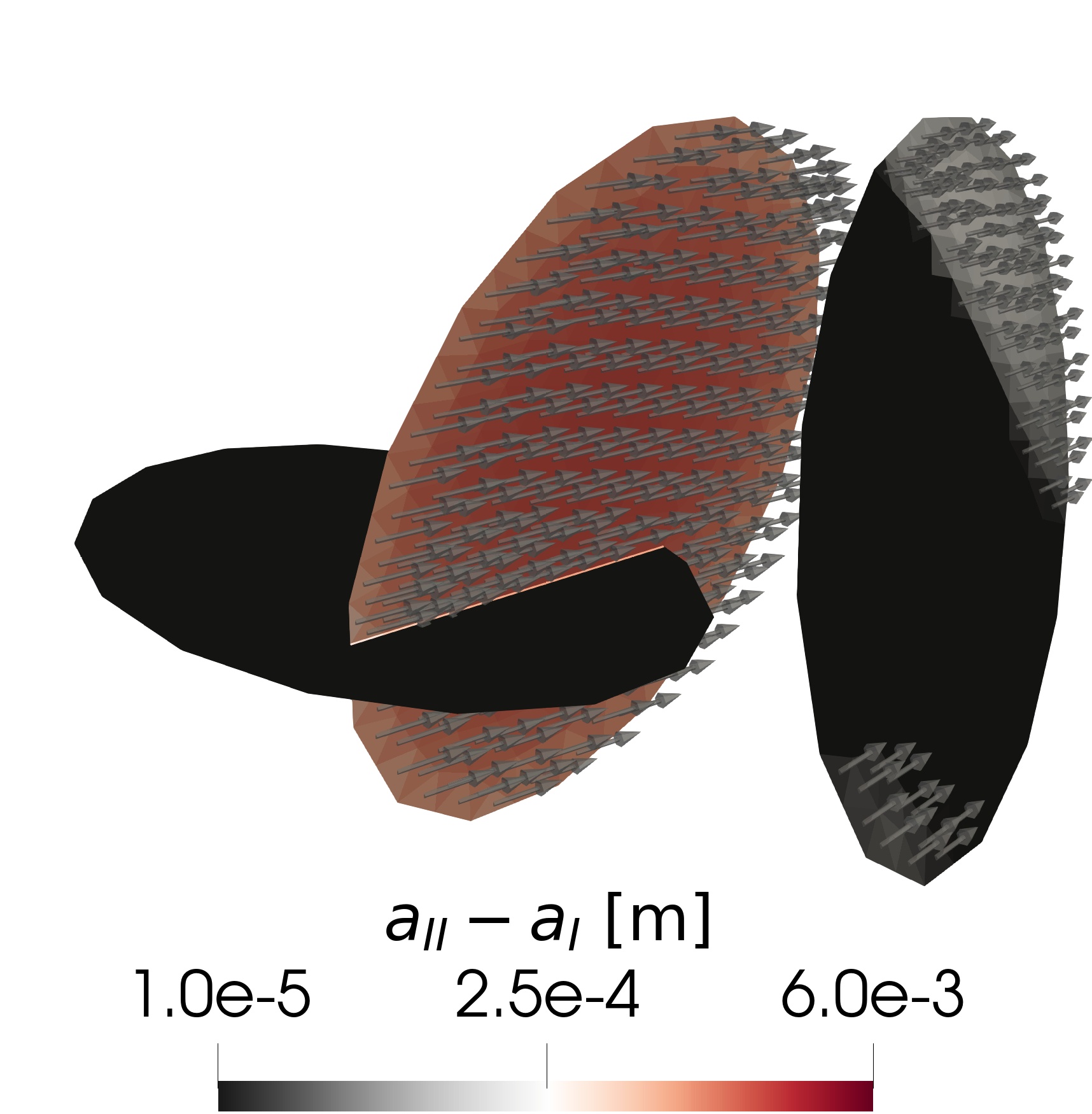}

\includegraphics[width=.42\textwidth]{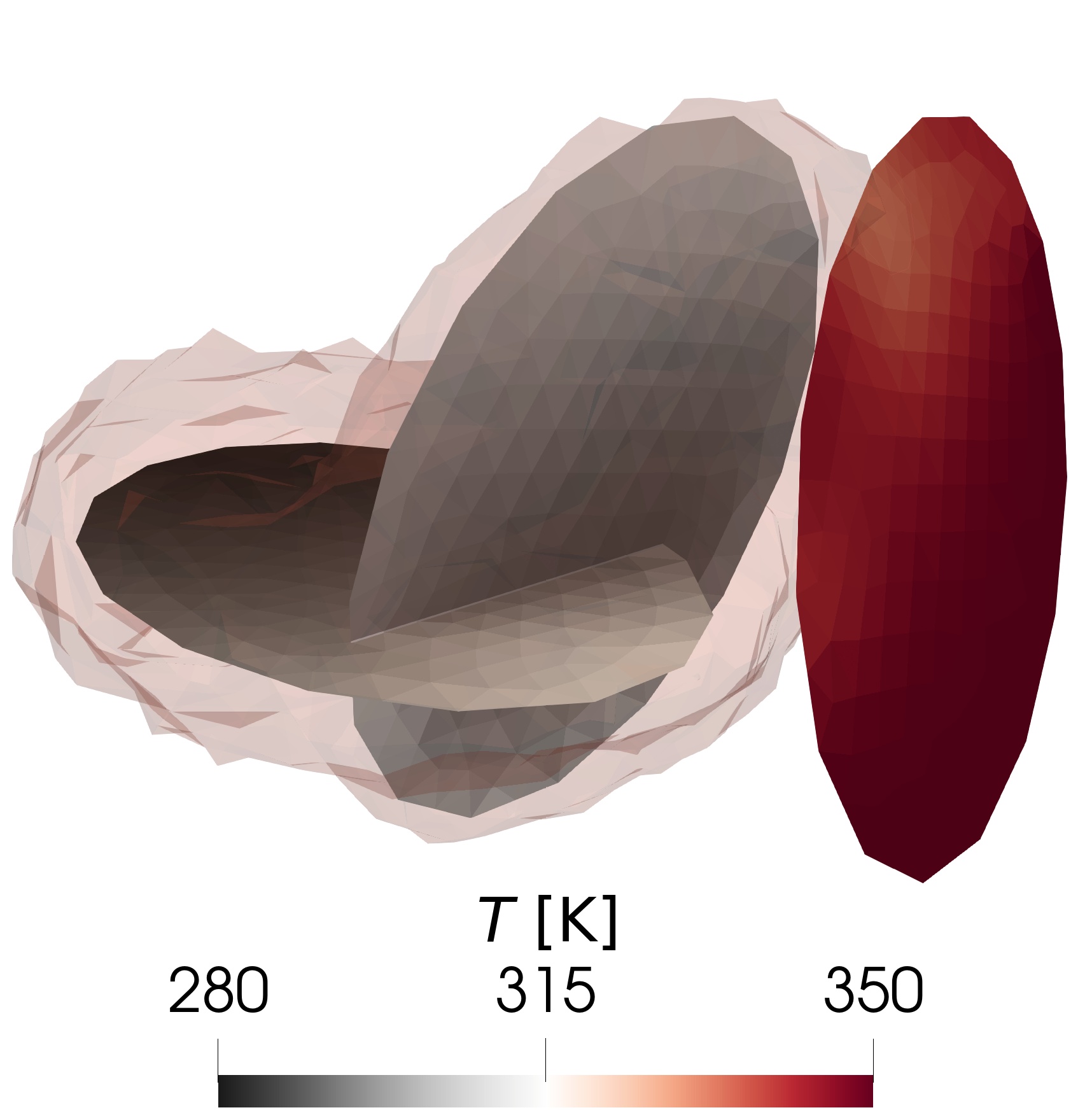}
\hfill
\includegraphics[width=.42\textwidth]{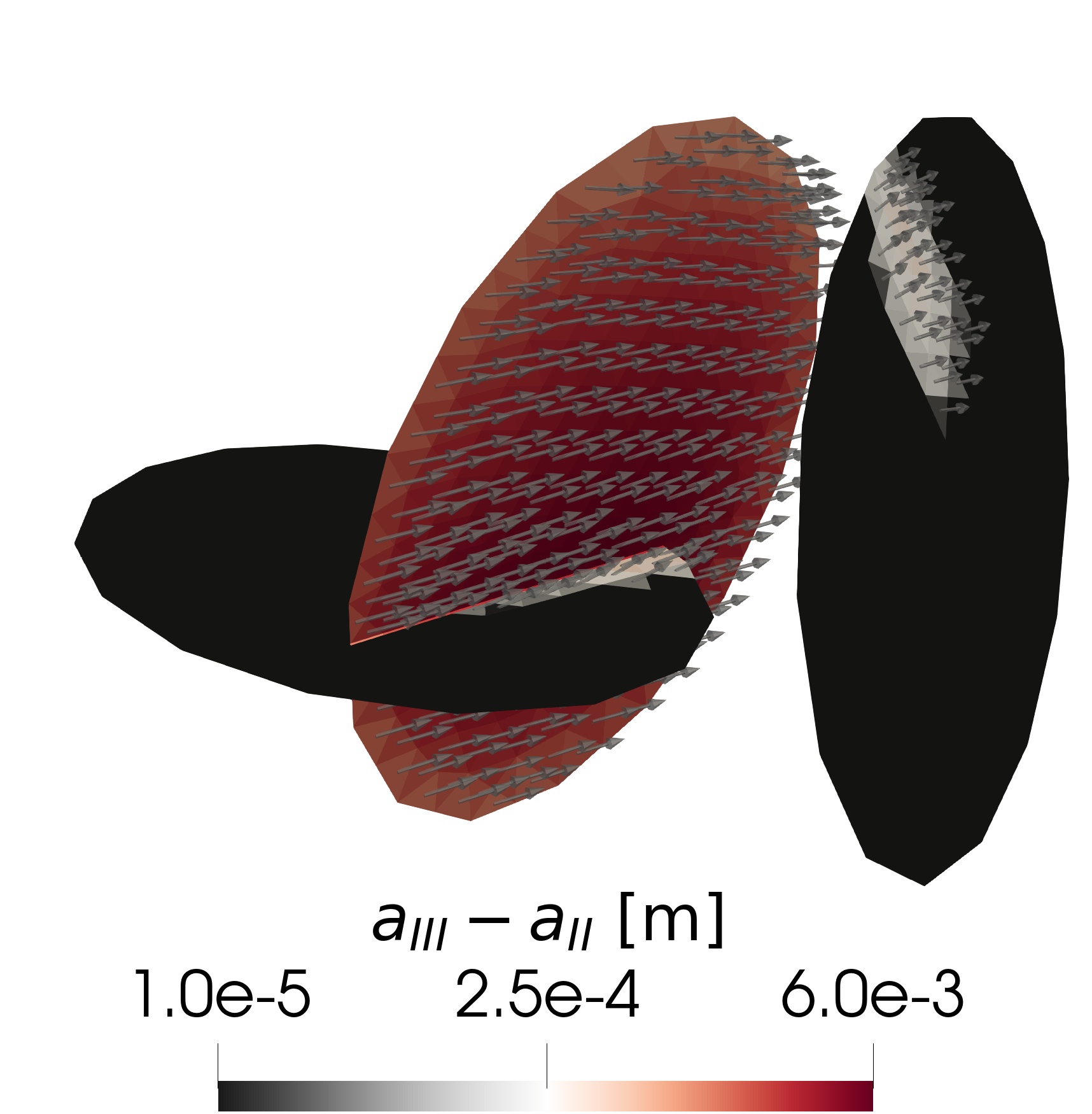}
\caption{Example 3: Left:  Deformation state according to the cumulative displacement jumps at the end of phase I (top). Perturbation from hydrostatic pressure at the end of phase II for the fractures and matrix cells satisfying $\pressure-\pressure[H]>\SI{5e5}{\pascal}$ (centre). Final fracture temperature and matrix cells satisfying $\temperature-\temperature[0]<\SI{-15}{\kelvin}$ (bottom).  
Right: Aperture increments and displacement jump direction for each of the three phases top to bottom.
The arrows show how the side of the fracture which is in view is displaced relative to the side not in view.
Note that the logarithmic scale for the aperture increments is truncated at \SI{1e-5}{m} for visualisation purposes. }
\label{fig:ex3_final_solutions}
\end{figure}

 Figure \ref{fig:ex3_final_solutions} shows spatial plots of pressure, temperature, aperture, displacement jumps and deformation state.
The plots demonstrate the model's cell-wise spatial resolution of the dynamics both in fractures and matrix.
For all three phases, displacement jumps are oriented in agreement with the background stress field and are very closely aligned. The only cells in \openSet are around the intersection in phase III. For the remaining cells, the (relatively small) normal components of the orientation arrows are due solely to shear dilation.

During phase II, aperture increments are most pronounced on fracture 2, which has no wells. However,  its intersection with fracture 1 where injection occurs leads to a significant pressure increase. Despite negligible pressure perturbation in fracture 3, some slip is observed due to stress redistribution following the deformation of fracture 2. The location of the slip in fracture 3, away from the stress shadow of fracture 2, highlights the complex mechanical interplay between fractures in a network. 

During phase III, some displacement jumps are induced in fracture 1 close to the intersection, whereas there is significant aperture increase throughout fracture 2 as a result of cooling of the surrounding rock. Along fracture 3, the deforming region is different from the previous phase, with displacement occurring in the region closest to fracture 2, where the surrounding matrix has been cooled the most. This conforms with the observations in Section \ref{sec:ex1} of aperture increases in regions of fluid entry or departure from the fractures.

\section{Conclusion}\label{sec:conclusion}
A model for fully coupled thermo-hydro-mechanical processes in porous media with deforming fractures is presented. 
Using the discrete--fracture--matrix approach, the matrix, the fractures and the fracture intersections are represented by subdomains of different dimensions connected by 
interfaces in a mixed-dimensional model.
Balance equations for energy and mass in all subdomains are coupled by fluxes on the interfaces, while the momentum balance in the matrix and traction balance and non-penetration for the fracture surfaces are coupled through interface displacements.  These governing equations are supplemented with constitutive laws, including a Coulomb type friction law and a linear shear dilation relation for the fractures. For the latter, a novel model consistently coupling slip and shear dilation of the fractures with the stress response of  the matrix is presented. 
The resulting set of model equations is discretised using cell-centred multi-point finite volume schemes and a semismooth Newton method for fracture deformation and solved fully coupled.

The model and its implementation are verified through a convergence study displaying first-order convergence for all primary variables and subdomains, except for the expected local reduction of convergence in the transition between contact regimes.
An exploration of three different shear-dilation models reveals significant discrepancies, demonstrating the importance of accurate and consistent modelling of the underlying physical mechanisms and their couplings.

Investigations of 2d and 3d examples identify a mechanism by which cooling-induced shear dilation preferentially occurs in regions where fluid leaves or enters a fracture.
The investigations also show the complexity of the process--structure interactions which may arise: in particular, how fracture deformation and resulting fracture dilation is induced by both mechanical, hydraulic and thermal driving forces. This confirms the need for models which explicitly incorporate  all relevant processes and structural features as well as the resulting process--structure interactions.
Furthermore, it demonstrates the proposed model's prowess in capturing such highly complex interactions and identifying their governing mechanisms.
Extensions such as chemical processes and more advanced friction models and dilation relations could readily be accommodated in the applied mixed-dimensional framework.

\section{Acknowledgements}
The authors thank two anonymous reviewers for comments and suggestions which contributed to improving the paper.

Funding: This work was supported by the Research Council of Norway and Equinor ASA through grant number 267908. \is{ I. Stefansson and I. Berre also acknowledge funding from the VISTA programme, The Norwegian Academy of Science  and Letters.}
\bibliographystyle{elsarticle-num-names}
\bibliography{bibliography}

\end{document}